\documentclass[12pt]{article}
\usepackage{amssymb,amsfonts }
\usepackage{ytableau}
\usepackage{amssymb}
\usepackage[margin=1in]{geometry}
\usepackage{mathrsfs}
\usepackage{stmaryrd}
\usepackage{amsfonts,amsmath,amssymb,amscd}
\usepackage{shadow}
\usepackage{graphicx}
\usepackage{color}
\usepackage{longtable}
\usepackage{pstricks,multido}
\usepackage{hyperref}
\usepackage{qtree}
\usepackage{tikz}\usetikzlibrary{matrix,arrows}
\usepackage[utf8]{inputenc}
\usepackage[T1]{fontenc}
\allowdisplaybreaks

\newcommand{\rmnum}[1]{\romannumeral #1}
\newtheorem{thm}{Theorem}[section]
\newtheorem{lem}[thm]{Lemma}

\newtheorem{example}{Example}[section]
\newtheorem{definition}{Definition}[section]

\newtheorem{conj}{Conjecture}[section]

\def\qed{\hfill \rule{4pt}{7pt}}

\def\pf{\noindent {\it{Proof.} \hskip 2pt}}

\numberwithin{equation}{section}

\begin{document}

\begin{center}
{\Large\bf On self-conjugate   $(s, s+1,\ldots, s+k)$-core partitions   }
\end{center}

\begin{center}
{\small Sherry H.F. Yan, Yao Yu, Hao Zhou}

Department of Mathematics\\
Zhejiang Normal University\\
 Jinhua 321004, P.R. China

 huifangyan@hotmail.com

\end{center}
\noindent {\bf Abstract.} Simultaneous core partitions have been widely studied since  Anderson's work  on the enumeration of $(s,t)$-core partitions. Amdeberhan and Leven    showed  that  the number of $(s,s+1, \ldots, s+k)$-core partitions is equal to the number of $(s, k)$-Dyck paths. In this paper, we prove that     self-conjugate  $(s,s+1,  \ldots, s+k)$-core partitions are equinumerous with symmetric  $(s, k)$-Dyck paths, confirming  a   conjecture  posed by Cho, Huh and Sohn.



%===========================================================================

\section{Introduction}
 A {\em partition} $\lambda$ of a positive integer $n$ is defined to
be a sequence of nonnegative  integers $(\lambda_1, \lambda_2, \ldots, \lambda_m)$ such that $\lambda_1+\lambda_2+\cdots+
\lambda_m=n$ and $\lambda_1\geq \lambda_2\cdots \geq \lambda_m$.
    The {\em Young diagram} of $\lambda$ is defined to be an up- and left- justified array of $n$ boxes with $\lambda_i$ boxes in the $i$-th row.  The {\em conjugate } of $\lambda$, denoted by $\lambda'$, is the partition whose Young diagram is the reflection along the main diagonal of $\lambda$'s diagram, and $\lambda$ is said to be {\em self-conjugate} if   $\lambda=\lambda'$. The {\em hook} of each box $B$ in $\lambda$ consists of the box $B$ itself and boxes directly to the right and directly below $B$. For example, Figure \ref{Young} illustrates
the Young diagram and the hook lengths of  a self-conjugate partition $\lambda=(5,3,3,1,1)$.

\begin{figure}[h!]
\begin{center}
 \begin{ytableau}
     9& 6 & 5 & 2&1\\
     6 & 3 & 2\\
    5 & 2 &1\\
     2\\
     1
 \end{ytableau}
 \end{center}
\caption{ The Young  diagram  of  a self-conjugate partition $\lambda=(5,3,3,1,1)$.}\label{Young}
\end{figure}

For a positive integer $s$, a partition is said to be an {\em $s$-core partition}, or simply an {\em $s$-core},  if it contains no box whose hook length is a multiple of $s$.  Furthermore,  a partition is said to be {\em $(t_1, t_2, \ldots, t_m)$-core} if it is simultaneously $t_1$-core, $t_2$-core, \ldots,  and $t_m$-core.  Anderson \cite{And}  showed that the number of $(s,t)$-core partitions is the
rational Catalan number  ${1\over s+t}{s+t\choose s}$
when $s$ and $t$ are coprime to each other.

  Simultaneous core partitions have attracted much attention since  Anderson's work  \cite{And}.  Various  results on the enumeration of such   partitions are archived in  \cite{Agg,Arm, Baek, Chen, Joh, Sta2, Str, Wang, Xiong, Yan, Zaleski,Zaleski2}.

  Ford, Mai and Sze \cite{Ford} characterized the set of hook lengths
of diagonal cells in self-conjugate  $(s,t)$-core partitions, and they showed that the number of
self-conjugate $(s,t)$-core partitions is given by ${\lfloor {s\over 2}\rfloor+\lfloor {t\over 2}\rfloor\choose \lfloor {s\over 2}\rfloor}$.
Concerning the enumeration of $(s,s+1, s+2)$-core partitions, Amdeberhan \cite{Amd}  conjectured that the number of $(s, s+1, s+2)$-core  partitions is equal to the number of Motzkin paths of length $s$. This conjecture was first  confirmed    by Yang, Zhong and Zhou \cite{Yang}.
Amdeberhan and Leven  \cite{Amd-lev} further proved that  the number of $(s,s+1, \ldots, s+k)$-core partitions is equal to the number of $(s, k)$-Dyck paths.
The  largest  size of such  partitions was derived by Xiong \cite{Xiong2}.

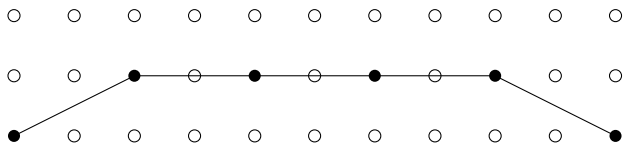
\begin{figure}[h!]
\begin{center}
\begin{picture}(250,100)
\setlength{\unitlength}{8mm}

 \put(0,0){\circle*{0.2}} \put(1,0){\circle{0.2}} \put(2,0){\circle{0.2}} \put(3,0){\circle{0.2}} \put(4,0){\circle{0.2}} \put(5,0){\circle{0.2}}
  \put(6,0){\circle{0.2}} \put(7,0){\circle{0.2}} \put(8,0){\circle{0.2}}  \put(9,0){\circle{0.2}} \put(10,0){\circle*{0.2}}

  \put(0,1){\circle{0.2}} \put(1,1){\circle{0.2}} \put(2,1){\circle*{0.2}} \put(3,1){\circle{0.2}} \put(4,1){\circle*{0.2}} \put(5,1){\circle{0.2}}
  \put(6,1){\circle*{0.2}} \put(7,1){\circle{0.2}} \put(8,1){\circle*{0.2}}  \put(9,1){\circle{0.2}} \put(10,1){\circle{0.2}}

\put(0,2){\circle{0.2}} \put(1,2){\circle{0.2}} \put(2,2){\circle{0.2}} \put(3,2){\circle{0.2}} \put(4,2){\circle{0.2}} \put(5,2){\circle{0.2}}
  \put(6,2){\circle{0.2}} \put(7,2){\circle{0.2}} \put(8,2){\circle{0.2}}  \put(9,2){\circle{0.2}} \put(10,2){\circle{0.2}}

\put(0,0){\line(2,1){2}}\put(2,1){\line(2,0){2}}
 \put(4,1){\line(2,0){2}}\put(6,1){\line(2,0){2}}
\put(8,1){\line(2,-1){2}}

\end{picture}
\end{center}
\caption{ A symmetric $(10,4)$-Dyck path.}\label{sym}
\end{figure}

A  { ballot $(s,k)$-path} of height $n$  is a lattice path  from $(0,0)$ to $(s,n)$
 using up steps $U_k=({k\over 2},1)$, down steps $D_k=({k\over 2},-1) $ and horizontal steps $H_{\ell}=(\ell,0)$ for some positive integer $\ell$ with  $1\leq \ell <k$,  and never
  lying below the $x$-axis.  A  ballot $(s,k)$-path of height $0$ is said to be an {\em (s,k)-Dyck} path. Notice that an (s,2)-Dyck  path is just a Motzkin path of length $s$. We say that an
$(s, k)$-Dyck path is {\em symmetric} if its reflection about the line $x={s\over 2}$ is itself. For example, a  symmetric $(10,4)$-Dyck path $P=U_2H_2H_2H_2D_2$ is shown in Figure \ref{sym}.

In this paper, we are mainly concerned with  the enumeration of  self-conjugate  $(s,s+1,  \ldots, s+k)$-core partitions, confirming the following  conjecture  posed by Cho, Huh and Sohn \cite{Cho}.

\begin{conj}\label{con1}
 For given positive integers $s$ and $k$, the number of self-conjugate $(s, s +
1, \ldots, s + k)$-core partitions is equal to the number of symmetric $(s, k)$-Dyck paths.
\end{conj}

 The case when $k=2$ has been confirmed by Cho, Huh and Sohn \cite{Cho}. To be more precise, they showed that the number of self-conjugate $(s, s +
1, s+2)$-core partitions is given by
$$
\sum_{i\geq 0}{\lfloor {s\over 2}\rfloor \choose i}{i\choose \lfloor {i\over 2}\rfloor  },
$$
which counts the    number of symmetric Motzkin paths of length $s$.

\section{Proof of Conjecture \ref{con1}}

In this section,  we aim to confirm Conjecture \ref{con1}. We begin with some definitions and notations.

  Following  the poset terminology given
by Stanley \cite{Sta},  for two elements $x$ and $y$ in a poset $(P, \prec)$, we say that $y$ {\em covers} $x$ if $x\prec y$ and there exists no element $z\in P$ such that $x\prec z\prec y$.
An element $x$ in $P$ is said to be of {\em rank} $s$ if it  covers an element of rank $s-1$. Note that the elements of rank $0$   in $P$ are just the minimal elements.
The {\em Hasse diagram} of a finite poset $P$ is a graph whose vertices are the elements of $P$,  whose edges are the cover relations, and such that if $y$  covers $x$ then there is an edge connecting $x$ and $y$ and
 $y$ is placed above $x$.
     An {\em order ideal} of $P$ is a subset $I$ such that if any $y\in I$ and $x\prec y$ in $P$, then $x\in I$.

For a partition $\lambda$, let $MD(\lambda)$ denote the set of      main diagonal hook lengths. It is apparent that  $MD(\lambda)$ is a set of distinct odd integers when $\lambda$ is self-conjugate. In the following theorem,  Ford, Mai and Sze \cite{Ford} provided a characterization of the set $MD(\lambda)$ of any self-conjugate $t$-core partition $\lambda$.
\begin{thm} { \upshape   (  \cite{Ford}) }\label{thford}
Let $\lambda$ be a self-conjugate partition. Then $\lambda$ is $t$-core if and only if  $MD(\lambda)$ verifies the following properties.
\begin{itemize}
\item[{\upshape (\rmnum{1})} ]  $MD(\lambda)$ is a set of  distinct odd integers;
\item[{\upshape (\rmnum{2})} ] If $h\in MD(\lambda)$ with $h>2t$, then $h-2t\in MD(\lambda)$;
\item[{\upshape (\rmnum{3})} ]  If $h_1, h_2\in MD(\lambda)$, then $  h_1+h_2$ is not a multiple of $2t$.
\end{itemize}
\end{thm}

 \begin{definition}
 Define  $L(s,k)=\bigcup_{j\geq 0} A_j,$ and $R(s,k)=\bigcup_{j\geq 0} B_j,$
  where
 $$
 A_j=\{2i-1+2sj\mid jk+1\leq i\leq \lfloor{s\over 2}\rfloor \},
 $$

 and
 $$
 B_j=\{2i-1+2sj\mid  jk+ \lceil{s+k\over 2} \rceil+1 \leq   i\leq s\}.
 $$
 \end{definition}

 \begin{lem}\label{lem3}
Let $x$ be a  positive odd  integer with $x>2s$. Then  $x$ is contained in $L(s,k) $ (resp. $R(s,k)$) if and only if  $x-2t\in L(s,k) $ (resp. $x-2t\in R(s,k) $ ) for all $s\leq t\leq s+k$.
\end{lem}
\pf
 The necessity follows directly from the definition of $L(s,k) $ (resp. $R(s,k)$).
Now we proceed to prove the sufficiency part of the theorem. Let   $x=2i-1+2js$  with  $1\leq i\leq s $ and $j\geq 1$. Assume that  $x-2t\in L(s,k) $ for all $s\leq t\leq s+k$. We proceed to show that $x\in L(s,k)$.  We claim that $i>k$. Since $x-2s-2k\in L(s,k)$, we have either $i>k$ or $j\geq 2$. Notice that   $x-2s=2i-1+2(j-1)s$. Then   $x-2s\in L(s,k)$ implies that $\lfloor{s\over 2}\rfloor\geq i\geq k(j-1)+1$. If $j\geq 2$, then  we also have $i>k$.      Then $x-2s-2k=2(i-k)-1+2(j-1)s\in L(s,k) $  implies that $i-k \geq k(j-1)+ 1$, which is equivalent to  $i\geq kj+1$. On the other hand,    $x-2s=2i-1+2(j-1)s\in L(s,k)$ implies that $  i  \leq \lfloor{s\over 2}\rfloor$.  Hence we have $x\in L(s,k) $.

    Assume that $x-2t\in R(s,k) $ for all $s\leq t\leq s+k$. We proceed to verify that $x\in R(s,k)$ by considering the following  two cases.

     If $j\geq 2$, then $x-2s=2i-1+2(j-1)s\in R(s,k)$ implies that $i-k(j-1)=i-kj+k\geq \lceil{s+k\over 2}\rceil+1$, which leads to   $  i- k \geq  \lceil{s+k\over 2}\rceil+1$.     Then $x-2s-2k=2(i-k)-1+2(j-1)s\in R(s,k) $  implies that $i-kj=i-k-k(j-1)\geq \lceil{s+k\over 2}\rceil+1$. Hence, we have $x\in R(s,k) $.

       If $j=1$, then $x-2s-2k=2i-1-2k\in R(s,k)$ implies that $i-k\geq \lceil{s+k\over 2}\rceil+1 $, which is equivalent to $i\geq k+ \lceil{s+k\over 2}\rceil+1$.  Thus,  we   have $x\in R(s,k)$.   This completes the proof. \qed

 \begin{definition}
Let $$P(s,k)=L(s,k)   \cup  R(s,k) $$
 with the following poset structure:  if $x, y\in P(s,k) $, then $y$ covers $x$ if and only if    $y=x+2s+2t$ for some $0\leq t\leq k$.
 \end{definition}
  Lemma \ref{lem3} tells us that $P(s,k)$ can be viewed as the disjoint union of the posets $L(s,k)$ and $R(s,k)$.   In the Hasse diagram of $L(s, k)$ (resp. $R(s, k)$),
 each element $x$ of rank $\ell$    covers exactly $ k+1$ elements   $ x-2s, x-2s-2$, $\ldots, x-2s-2k $ of rank $\ell-1$ for all $\ell\geq 1$. For example,  the Hasse   diagrams of $P(20,3)$ and $ P(20,4)$ are illustrated in Figure \ref{P{20,3}} and \ref{P{20,4}}, respectively.

An order ideal $I$ of $P(s,k)$ is said to be  {\em nice } if  $I$ does not contain  two elements $h_1, h_2$ such that  $h_1+h_2\in  \{2s, 2s+2, \ldots, 2s+2k\}$.

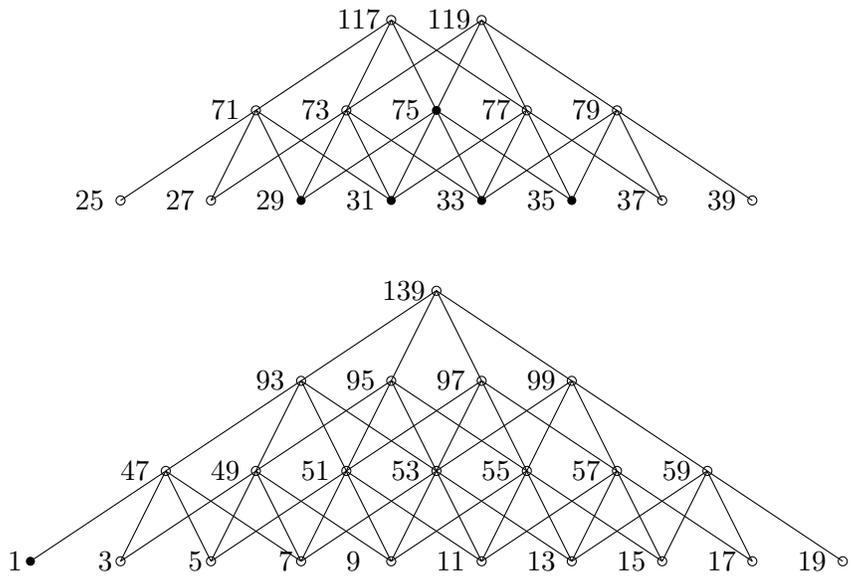
\begin{figure}[h!]
\begin{center}
\begin{picture}(230,150)
\setlength{\unitlength}{6mm}

 \put(0,0){\circle*{0.2}}
 \put(2,0){\circle{0.2}}
 \put(4,0){\circle{0.2}}
 \put(6,0){\circle{0.2}}
 \put(8,0){\circle{0.2}}
 \put(10,0){\circle{0.2}}
  \put(12,0){\circle{0.2}}
 \put(14,0){\circle{0.2}}
\put(16,0){\circle{0.2}}
\put(18,0){\circle{0.2}}

\put(3,2){\circle{0.2}}\put(5,2){\circle{0.2}}\put(7,2){\circle{0.2}}
\put(9,2){\circle{0.2}}\put(11,2){\circle{0.2}}\put(13,2){\circle{0.2}}
\put(15,2){\circle{0.2}}

\put(6,4){\circle{0.2}}
\put(8,4){\circle{0.2}}\put(10,4){\circle{0.2}}
\put(12,4){\circle{0.2}}\put(9,6){\circle{0.2}}

\put(0,0){\line(3,2){3}}\put(2,0){\line(3,2){3}}\put(4,0){\line(3,2){3}}
\put(6,0){\line(3,2){3}}\put(8,0){\line(3,2){3}}\put(10,0){\line(3,2){3}}
\put(12,0){\line(3,2){3}}
\put(3,2){\line(3,2){3}}\put(5,2){\line(3,2){3}}\put(7,2){\line(3,2){3}}
\put(9,2){\line(3,2){3}}\put(6,4){\line(3,2){3}}
\put(2,8){\line(3,2){3}}\put(4,8){\line(3,2){3}}\put(6,8){\line(3,2){3}}
\put(8,8){\line(3,2){3}}\put(10,8){\line(3,2){3}}
\put(5,10){\line(3,2){3}}\put(7,10){\line(3,2){3}}

\put(2,0){\line(1,2){1}}\put(4,0){\line(1,2){1}}\put(6,0){\line(1,2){1}}
\put(8,0){\line(1,2){1}}\put(10,0){\line(1,2){1}}\put(12,0){\line(1,2){1}}
\put(14,0){\line(1,2){1}}
\put(5,2){\line(1,2){1}}\put(7,2){\line(1,2){1}}\put(9,2){\line(1,2){1}}
\put(11,2){\line(1,2){1}}\put(8,4){\line(1,2){1}}
\put(4,8){\line(1,2){1}}\put(6,8){\line(1,2){1}}\put(8,8){\line(1,2){1}}
\put(10,8){\line(1,2){1}}\put(12,8){\line(1,2){1}}\put(7,10){\line(1,2){1}}
\put(9,10){\line(1,2){1}}

\put(4,0){\line(-1,2){1}}\put(6,0){\line(-1,2){1}}\put(8,0){\line(-1,2){1}}
\put(10,0){\line(-1,2){1}}\put(12,0){\line(-1,2){1}}\put(14,0){\line(-1,2){1}}
\put(16,0){\line(-1,2){1}}
\put(7,2){\line(-1,2){1}}\put(9,2){\line(-1,2){1}}\put(11,2){\line(-1,2){1}}
\put(13,2){\line(-1,2){1}}\put(10,4){\line(-1,2){1}}
\put(6,8){\line(-1,2){1}}\put(8,8){\line(-1,2){1}}\put(10,8){\line(-1,2){1}}
\put(12,8){\line(-1,2){1}}\put(14,8){\line(-1,2){1}}
\put(9,10){\line(-1,2){1}}\put(11,10){\line(-1,2){1}}

\put(6,0){\line(-3,2){3}}\put(8,0){\line(-3,2){3}}\put(10,0){\line(-3,2){3}}
\put(12,0){\line(-3,2){3}}\put(14,0){\line(-3,2){3}}\put(16,0){\line(-3,2){3}}
\put(18,0){\line(-3,2){3}}
\put(9,2){\line(-3,2){3}}\put(11,2){\line(-3,2){3}}\put(13,2){\line(-3,2){3}}
\put(15,2){\line(-3,2){3}}\put(12,4){\line(-3,2){3}}
\put(8,8){\line(-3,2){3}}\put(10,8){\line(-3,2){3}}\put(12,8){\line(-3,2){3}}
\put(14,8){\line(-3,2){3}}\put(16,8){\line(-3,2){3}}
\put(11,10){\line(-3,2){3}}\put(13,10){\line(-3,2){3}}

\put(2,8){\circle{0.2}}\put(4,8){\circle{0.2}}\put(6,8){\circle*{0.2}}\put(8,8){\circle*{0.2}}
\put(10,8){\circle*{0.2}}\put(12,8){\circle*{0.2}}\put(14,8){\circle{0.2}}\put(16,8){\circle{0.2}}

\put(5,10){\circle{0.2}}
\put(7,10){\circle{0.2}}\put(9,10){\circle*{0.2}}\put(11,10){\circle{0.2}}\put(13,10){\circle{0.2}}
\put(8,12){\circle{0.2}}\put(10,12){\circle{0.2}}

 \put(-0.5, -0.2){\small$1$}\put(1.5, -0.2){\small$3$}\put(3.5, -0.2){\small$5$}\put(5.5, -0.2){\small$7$}\put(7, -0.2){\small$9$}
 \put(9, -0.2){\small$11$} \put(11, -0.2){\small$13$} \put(13, -0.2){\small$15$} \put(15, -0.2){\small$17$} \put(17, -0.2){\small$19$}

 \put(2, 1.8){\small$47$}\put(4, 1.8){\small$49$}\put(6, 1.8){\small$51$}
 \put(8, 1.8){\small$53$}\put(10, 1.8){\small$55$}\put(12, 1.8){\small$57$}
 \put(14, 1.8){\small$59$}
 \put(5, 3.8){\small$93$}\put(7, 3.8){\small$95$}\put(9, 3.8){\small$97$}
 \put(11, 3.8){\small$99$}\put(7.8, 5.8){\small$139$}

\put(1,7.8){\small$25$}\put(3,7.8){\small$27$}\put(5,7.8){\small$29$}
\put(7,7.8){\small$31$}\put(9,7.8){\small$33$}\put(11,7.8){\small$35$}
\put(13,7.8){\small$37$}\put(15,7.8){\small$39$}
\put(4,9.8){\small$71$}\put(6,9.8){\small$73$}\put(8,9.8){\small$75$}
\put(10,9.8){\small$77$}\put(12,9.8){\small$79$}
\put(6.8,11.8){\small$117$}\put(8.8,11.8){\small$119$}

\end{picture}
\end{center}
\caption{ The Hasse diagram of the poset  $P(20,3)$.}\label{P{20,3}}
\end{figure}

 \begin{figure}[h!]
\begin{center}
\begin{picture}(230,150)
\setlength{\unitlength}{6mm}

 \put(0,0){\circle*{0.2}}
 \put(2,0){\circle*{0.2}}
 \put(4,0){\circle*{0.2}}
 \put(6,0){\circle*{0.2}}
 \put(8,0){\circle*{0.2}}
 \put(10,0){\circle{0.2}}
  \put(12,0){\circle{0.2}}
 \put(14,0){\circle{0.2}}
\put(16,0){\circle{0.2}}
\put(18,0){\circle{0.2}}

\put(4,2){\circle*{0.2}}\put(6,2){\circle{0.2}}\put(8,2){\circle{0.2}}
\put(10,2){\circle{0.2}}\put(12,2){\circle{0.2}}\put(14,2){\circle{0.2}}

\put(8,4){\circle{0.2}}\put(10,4){\circle{0.2}}

\put(8,0){\line(-4,2){4}}
\put(10,0){\line(-4,2){4}}\put(12,0){\line(-4,2){4}}\put(14,0){\line(-4,2){4}}
\put(16,0){\line(-4,2){4}}\put(18,0){\line(-4,2){4}}

\put(12,2){\line(-4,2){4}}
\put(14,2){\line(-4,2){4}}

\put(6,0){\line(-2,2){2}}\put(8,0){\line(-2,2){2}}\put(10,0){\line(-2,2){2}}
\put(12,0){\line(-2,2){2}}\put(14,0){\line(-2,2){2}}\put(16,0){\line(-2,2){2}}

\put(10,2){\line(-2,2){2}}\put(12,2){\line(-2,2){2}}\put(14,2)

 \put(0,0){\line(4,2){4}} \put(2,0){\line(4,2){4}} \put(4,0){\line(4,2){4}} \put(6,0){\line(4,2){4}} \put(8,0){\line(4,2){4}} \put(10,0){\line(4,2){4}}
 \put(4,2){\line(4,2){4}}\put(6,2){\line(4,2){4}}

  \put(2,0){\line(2,2){2}}
   \put(4,0){\line(2,2){2}}
 \put(6,0){\line(2,2){2}}
 \put(8,0){\line(2,2){2}}
 \put(10,0){\line(2,2){2}}
 \put(12,0){\line(2,2){2}}

\put(6,2){\line(2,2){2}}\put(8,2){\line(2,2){2}}

 \put(4,0){\line(0,2){2}} \put(6,0){\line(0,2){2}} \put(8,0){\line(0,2){2}} \put(10,0){\line(0,2){2}} \put(12,0){\line(0,2){2}} \put(14,0){\line(0,2){2}}

\put(8,2){\line(0,2){2}} \put(10,2){\line(0,2){2}}

 \put(-0.5, -0.2){\small$1$}\put(1.5, -0.2){\small$3$}\put(3.5, -0.2){\small$5$}\put(5.5, -0.2){\small$7$}\put(7, -0.2){\small$9$}
 \put(9, -0.2){\small$11$} \put(11, -0.2){\small$13$} \put(13, -0.2){\small$15$} \put(15, -0.2){\small$17$} \put(17, -0.2){\small$19$}

  \put(3,1.8){\small$49$}\put(5,1.8){\small$51$} \put(7,1.8){\small$53$}\put(9,1.8){\small$55$}\put(11,1.8){\small$57$}\put(13,1.8){\small$59$}

  \put(7, 3.8){\small$97$}\put(9, 3.8){\small$99$}

\put(2,6){\circle*{0.2}}\put(4,6){\circle{0.2}}\put(6,6){\circle*{0.2}}\put(8,6){\circle{0.2}}
\put(10,6){\circle{0.2}}\put(12,6){\circle{0.2}}\put(14,6){\circle{0.2}}\put(16,6){\circle{0.2}}

\put(6,8){\circle{0.2}}\put(8,8){\circle{0.2}}\put(10,8){\circle{0.2}}\put(12,8){\circle{0.2}}

\put(10,6){\line(-4,2){4}}\put(12,6){\line(-4,2){4}}\put(14,6){\line(-4,2){4}}
\put(16,6){\line(-4,2){4}}

\put(8,6){\line(-2,2){2}}\put(10,6){\line(-2,2){2}}\put(12,6){\line(-2,2){2}}
\put(14,6){\line(-2,2){2}}

\put(6,6){\line(0,2){2}}\put(8,6){\line(0,2){2}}\put(10,6){\line(0,2){2}}
\put(12,6){\line(0,2){2}}

\put(4,6){\line(2,2){2}}\put(6,6){\line(2,2){2}}
\put(8,6){\line(2,2){2}}\put(10,6){\line(2,2){2}}

\put(2,6){\line(4,2){4}}\put(4,6){\line(4,2){4}}\put(6,6){\line(4,2){4}}\put(8,6){\line(4,2){4}}
\put(1,5.8){\small$25$}\put(3,5.8){\small$27$}\put(5,5.8){\small$29$}
\put(7,5.8){\small$31$}\put(9,5.8){\small$33$}\put(11,5.8){\small$35$}
\put(13,5.8){\small$37$}\put(15,5.8){\small$39$}

\put(5,7.8){\small$73$}\put(7,7.8){\small$75$}\put(9,7.8){\small$77$}
\put(11,7.8){\small$79$}

\end{picture}
\end{center}
\caption{ The Hasse diagram of the poset  $P(20,4)$.}\label{P{20,4}}
\end{figure}

\begin{figure}[h!]
\begin{center}
\begin{picture}(350,220)
\setlength{\unitlength}{8mm}

 \put(0,0){\circle*{0.2}}
 \put(2,0){\circle{0.2}}
 \put(4,0){\circle{0.2}}
 \put(6,0){\circle{0.2}}
 \put(8,0){\circle{0.2}}
 \put(10,0){\circle{0.2}}
  \put(12,0){\circle{0.2}}
 \put(14,0){\circle{0.2}}
\put(16,0){\circle{0.2}}
\put(18,0){\circle{0.2}}

\put(3,2){\circle{0.2}}\put(5,2){\circle{0.2}}\put(7,2){\circle{0.2}}
\put(9,2){\circle{0.2}}\put(11,2){\circle{0.2}}\put(13,2){\circle{0.2}}
\put(15,2){\circle{0.2}}

\put(6,4){\circle{0.2}}
\put(8,4){\circle{0.2}}\put(10,4){\circle{0.2}}
\put(12,4){\circle{0.2}}\put(9,6){\circle{0.2}}

\put(0,0){\line(3,2){3}}\put(2,0){\line(3,2){3}}\put(4,0){\line(3,2){3}}
\put(6,0){\line(3,2){3}}\put(8,0){\line(3,2){3}}\put(10,0){\line(3,2){3}}
\put(12,0){\line(3,2){3}}
\put(3,2){\line(3,2){3}}\put(5,2){\line(3,2){3}}\put(7,2){\line(3,2){3}}
\put(9,2){\line(3,2){3}}\put(6,4){\line(3,2){3}}
\put(2,8){\line(3,2){3}}\put(4,8){\line(3,2){3}}\put(6,8){\line(3,2){3}}
\put(8,8){\line(3,2){3}}\put(10,8){\line(3,2){3}}
\put(5,10){\line(3,2){3}}\put(7,10){\line(3,2){3}}

\put(2,0){\line(1,2){1}}\put(4,0){\line(1,2){1}}\put(6,0){\line(1,2){1}}
\put(8,0){\line(1,2){1}}\put(10,0){\line(1,2){1}}\put(12,0){\line(1,2){1}}
\put(14,0){\line(1,2){1}}
\put(5,2){\line(1,2){1}}\put(7,2){\line(1,2){1}}\put(9,2){\line(1,2){1}}
\put(11,2){\line(1,2){1}}\put(8,4){\line(1,2){1}}
\put(4,8){\line(1,2){1}}\put(6,8){\line(1,2){1}}\put(8,8){\line(1,2){1}}
\put(10,8){\line(1,2){1}}\put(12,8){\line(1,2){1}}\put(7,10){\line(1,2){1}}
\put(9,10){\line(1,2){1}}

\put(4,0){\line(-1,2){1}}\put(6,0){\line(-1,2){1}}\put(8,0){\line(-1,2){1}}
\put(10,0){\line(-1,2){1}}\put(12,0){\line(-1,2){1}}\put(14,0){\line(-1,2){1}}
\put(16,0){\line(-1,2){1}}
\put(7,2){\line(-1,2){1}}\put(9,2){\line(-1,2){1}}\put(11,2){\line(-1,2){1}}
\put(13,2){\line(-1,2){1}}\put(10,4){\line(-1,2){1}}
\put(6,8){\line(-1,2){1}}\put(8,8){\line(-1,2){1}}\put(10,8){\line(-1,2){1}}
\put(12,8){\line(-1,2){1}}\put(14,8){\line(-1,2){1}}
\put(9,10){\line(-1,2){1}}\put(11,10){\line(-1,2){1}}

\put(6,0){\line(-3,2){3}}\put(8,0){\line(-3,2){3}}\put(10,0){\line(-3,2){3}}
\put(12,0){\line(-3,2){3}}\put(14,0){\line(-3,2){3}}\put(16,0){\line(-3,2){3}}
\put(18,0){\line(-3,2){3}}
\put(9,2){\line(-3,2){3}}\put(11,2){\line(-3,2){3}}\put(13,2){\line(-3,2){3}}
\put(15,2){\line(-3,2){3}}\put(12,4){\line(-3,2){3}}
\put(8,8){\line(-3,2){3}}\put(10,8){\line(-3,2){3}}\put(12,8){\line(-3,2){3}}
\put(14,8){\line(-3,2){3}}\put(16,8){\line(-3,2){3}}
\put(11,10){\line(-3,2){3}}\put(13,10){\line(-3,2){3}}

\put(2,8){\circle{0.2}}\put(4,8){\circle{0.2}}\put(6,8){\circle*{0.2}}\put(8,8){\circle*{0.2}}
\put(10,8){\circle*{0.2}}\put(12,8){\circle*{0.2}}\put(14,8){\circle{0.2}}\put(16,8){\circle{0.2}}

\put(5,10){\circle{0.2}}
\put(7,10){\circle{0.2}}\put(9,10){\circle*{0.2}}\put(11,10){\circle{0.2}}\put(13,10){\circle{0.2}}
\put(8,12){\circle{0.2}}\put(10,12){\circle{0.2}}

 \put(-1.2, -0.2){\small$(1,0)$}\put(0.8, -0.2){\small$(2,0)$}\put(2.8, -0.2){\small$(3,0)$}\put(4.8, -0.2){\small$(4,0)$}\put(6.8, -0.2){\small$(5,0)$}
 \put(8.8, -0.2){\small$(6,0)$} \put(10.8, -0.2){\small$(7,0)$} \put(12.8, -0.2){\small$(8,0)$} \put(14.8, -0.2){\small$(9,0)$} \put(16.5, -0.2){\small$(10,0)$}

 \put(1.4, 1.8){\small$(2.5,1)$}\put(3.4, 1.8){\small$(3.5,1)$}\put(5.4, 1.8){\small$(4.5,1)$}
 \put(7.4, 1.8){\small$(5.5,1)$}\put(9.4, 1.8){\small$(6.5,1)$}\put(11.4, 1.8){\small$(7.5,1)$}
 \put(13.4, 1.8){\small$(8.5,1)$}
 \put(4.7, 3.8){\small$(4,2)$}\put(6.5, 3.8){\small$(5,2)$}\put(8.5, 3.8){\small$(6,2)$}
 \put(10.5, 3.8){\small$(7,2)$}\put(7.2, 5.8){\small$(5.5,3)$}

\put(0.1,7.8){\small$(9.5,-1)$}\put(2.1,7.8){\small$(8.5,-1)$}\put(4.1,7.8){\small$(7.5,-1)$}
\put(6.1,7.8){\small$(6.5,-1)$}\put(8.1,7.8){\small$(5.5,-1)$}\put(10.1,7.8){\small$(4.5,-1)$}
\put(12.1,7.8){\small$(3.5,-1)$}\put(14.1,7.8){\small$(2.5,-1)$}
\put(3.3,9.8){\small$(8,-2)$}\put(5.3,9.8){\small$(7,-2)$}\put(7.3,9.8){\small$(6,-2)$}
\put(9.3,9.8){\small$(5,-2)$}\put(11.3,9.8){\small$(4,-2)$}
\put(7,12.2){\small$(6.5,-3)$}\put(9,12.2){\small$(5.5,-3)$}

\end{picture}
\end{center}
\caption{ The Hasse diagram of the poset  $P'(20,3)$.}\label{P'{20,3}}
\end{figure}

 \begin{figure}[h!]
\begin{center}
\begin{picture}(350,200)
\setlength{\unitlength}{8mm}

 \put(0,0){\circle*{0.2}}
 \put(2,0){\circle*{0.2}}
 \put(4,0){\circle*{0.2}}
 \put(6,0){\circle*{0.2}}
 \put(8,0){\circle*{0.2}}
 \put(10,0){\circle{0.2}}
  \put(12,0){\circle{0.2}}
 \put(14,0){\circle{0.2}}
\put(16,0){\circle{0.2}}
\put(18,0){\circle{0.2}}

\put(4,2){\circle*{0.2}}\put(6,2){\circle{0.2}}\put(8,2){\circle{0.2}}
\put(10,2){\circle{0.2}}\put(12,2){\circle{0.2}}\put(14,2){\circle{0.2}}

\put(8,4){\circle{0.2}}\put(10,4){\circle{0.2}}

\put(8,0){\line(-4,2){4}}
\put(10,0){\line(-4,2){4}}\put(12,0){\line(-4,2){4}}\put(14,0){\line(-4,2){4}}
\put(16,0){\line(-4,2){4}}\put(18,0){\line(-4,2){4}}

\put(12,2){\line(-4,2){4}}
\put(14,2){\line(-4,2){4}}

\put(6,0){\line(-2,2){2}}\put(8,0){\line(-2,2){2}}\put(10,0){\line(-2,2){2}}
\put(12,0){\line(-2,2){2}}\put(14,0){\line(-2,2){2}}\put(16,0){\line(-2,2){2}}

\put(10,2){\line(-2,2){2}}\put(12,2){\line(-2,2){2}}\put(14,2)

 \put(0,0){\line(4,2){4}} \put(2,0){\line(4,2){4}} \put(4,0){\line(4,2){4}} \put(6,0){\line(4,2){4}} \put(8,0){\line(4,2){4}} \put(10,0){\line(4,2){4}}
 \put(4,2){\line(4,2){4}}\put(6,2){\line(4,2){4}}

  \put(2,0){\line(2,2){2}}
   \put(4,0){\line(2,2){2}}
 \put(6,0){\line(2,2){2}}
 \put(8,0){\line(2,2){2}}
 \put(10,0){\line(2,2){2}}
 \put(12,0){\line(2,2){2}}

\put(6,2){\line(2,2){2}}\put(8,2){\line(2,2){2}}

 \put(4,0){\line(0,2){2}} \put(6,0){\line(0,2){2}} \put(8,0){\line(0,2){2}} \put(10,0){\line(0,2){2}} \put(12,0){\line(0,2){2}} \put(14,0){\line(0,2){2}}

\put(8,2){\line(0,2){2}} \put(10,2){\line(0,2){2}}

 \put(-1.5, -0.4){\small$(1,0)$}\put(0.5, -0.4){\small$(2,0)$}\put(2.5, -0.4){\small$(3,0)$}
 \put(4.5, -0.4){\small$(4,0)$}\put(6.5, -0.4){\small$(5,0)$}\put(8.5, -0.4){\small$(6,0)$}
 \put(10.5, -0.4){\small$(7,0)$}\put(12.5, -0.4){\small$(8,0)$}
 \put(14.5, -0.4){\small$(9,0)$}
\put(16.3, -0.4){\small$(10,0)$}

\put(2.5, 2){\small$(3,1)$}\put(4.5, 2){\small$(4,1)$}
\put(6.5, 2){\small$(5,1)$}\put(8.5, 2){\small$(6,1)$}\put(10.5, 2){\small$(7,1)$}
\put(12.5, 2){\small$(8,1)$}
\put(6.5, 4){\small$(5,2)$}\put(8.5, 4){\small$(6,2)$}

\put(2,6){\circle*{0.2}}\put(4,6){\circle{0.2}}\put(6,6){\circle*{0.2}}\put(8,6){\circle{0.2}}
\put(10,6){\circle{0.2}}\put(12,6){\circle{0.2}}\put(14,6){\circle{0.2}}\put(16,6){\circle{0.2}}

\put(6,8){\circle{0.2}}\put(8,8){\circle{0.2}}\put(10,8){\circle{0.2}}\put(12,8){\circle{0.2}}

\put(10,6){\line(-4,2){4}}\put(12,6){\line(-4,2){4}}\put(14,6){\line(-4,2){4}}
\put(16,6){\line(-4,2){4}}

\put(8,6){\line(-2,2){2}}\put(10,6){\line(-2,2){2}}\put(12,6){\line(-2,2){2}}
\put(14,6){\line(-2,2){2}}

\put(6,6){\line(0,2){2}}\put(8,6){\line(0,2){2}}\put(10,6){\line(0,2){2}}
\put(12,6){\line(0,2){2}}

\put(4,6){\line(2,2){2}}\put(6,6){\line(2,2){2}}
\put(8,6){\line(2,2){2}}\put(10,6){\line(2,2){2}}

\put(2,6){\line(4,2){4}}\put(4,6){\line(4,2){4}}\put(6,6){\line(4,2){4}}\put(8,6){\line(4,2){4}}

\put(0.5, 5.6){\small$(10,-1)$}\put(2.5, 5.6){\small$(9,-1)$}\put(4.5, 5.6){\small$(8,-1)$}
\put(6.5, 5.6){\small$(7,-1)$}\put(8.5, 5.6){\small$(6,-1)$}\put(10.5, 5.6){\small$(5,-1)$}
\put(12.5, 5.6){\small$(4,-1)$}\put(14.5, 5.6){\small$(3,-1)$}
\put(4.5, 8){\small$(8,-2)$}\put(6.5, 8){\small$(7,-2)$}\put(8.5, 8){\small$(6,-2)$}
\put(10.5, 8){\small$(5,-2)$}

\end{picture}
\end{center}
\caption{ The Hasse diagram of the poset  $P'(20,4)$.}\label{P'{20,4}}
\end{figure}

In the following theorem, we establish a correspondence between self-conjugate  $(s, s+1, \ldots, s+k)$-core partitions and the  nice order ideals of the poset $P(s,k) $.
\begin{thm}\label{main1}
Let $\lambda$ be a self-conjugate partition. Then $\lambda$ is    $(s, s+1, \ldots, s+k)$-core if and only if $MD(\lambda)$ is a nice order ideal of the poset $P(s,k) $.
 \end{thm}
\pf Assume that $\lambda$ is a self-conjugate  $(s, s+1, \ldots, s+k)$-core partition, we proceed to show that $MD(\lambda)$ is a nice order ideal of  $P(s,k)$. By   {\upshape (\rmnum{2})} and {\upshape (\rmnum{3})} of Theorem \ref{thford}, it suffices to show that    $MD(\lambda)$ is a subset of $P(s,k)$. Suppose to the contrary,  there exists an element $x$ of  $MD(\lambda)$ such that $x$ is not  contained in $P(s,k)$. Choose $x$ to the smallest such element. Examining  the definition of the poset $P(s,k)$, we have $x>2s$.  We claim that $x>2s+2k$. If not, suppose that $x=2s+2t+1$ for some $0\leq t\leq k-1$.  Then by {\upshape (\rmnum{2})} of Theorem \ref{thford}, $MD(\lambda)$ contains the element $y=1$. Then there are  two elements  $x=2s+2t+1$ and $y=1$ in $MD(\lambda)$ such  that $x+y=2s+2t+2$, contradicting {\upshape (\rmnum{3})}.   By {\upshape (\rmnum{2})}  of Theorem \ref{thford} and  the selection  of $x$, we have $x-2t\in P(s,k) $ for all $s\leq t\leq s+k$. From Lemma \ref{lem3}, it follows that $x\in P(s,k) $, which contradicts the selection of $x$. Hence we have reached the conclusion  that $MD(\lambda)$ is a subset of  $P(s,k)$.

Conversely, given a nice order ideal $I$ of $P(s,k)$, we wish to show that $I=MD(\lambda)$ for some self-conjugate $(s, s+1, \ldots, s+k)$-core partition $\lambda$. It suffices to prove that $I$ verifies properties {\upshape (\rmnum{1})}- {\upshape (\rmnum{3})} of Theorem \ref{thford}. It is easy to check that $I$ verifies properties {\upshape (\rmnum{1})} and   {\upshape (\rmnum{2})}. We proceed to show that $I$ also has property {\upshape (\rmnum{3})}.  Let    $h_1, h_2\in I$.  For any $1\leq t\leq k $, we have  $h_1=2p-1+2(s+t)q$ and $h_2=2p'-1+2(s+t)q'$ where $1\leq p,p'\leq s+t$, and $q,q'\geq 0$. Let $x=2p-1$ and $y=2p'-1$.  Since $I$ is a nice order ideal of $P(s,k)$, the elements $x$ and $y$ are contained in $I$ and $x+y\neq 2(s+t)$. This implies that $h_1+h_2$ is not a multiple of $2(s+t)$ for all $0\leq t\leq k$. Thus we have concluded that $I$ verifies  {\upshape (\rmnum{3})}.  This completes the proof. \qed

 \begin{definition}
 Define  $L'(s,k)=\bigcup_{j\geq 0} A'_j,$ and $R'(s,k)=\bigcup_{j\geq 0} B'_j,$
  where
 $$
 A'_j=\{(i-{k\over 2}j,j)\mid  jk+1\leq i\leq \lfloor{s\over 2}\rfloor,\,\, i\in \mathcal{Z}\},
 $$

 and
 $$
 B'_j=\{(s+1-i+{k\over 2}(j+1),-j-1)\mid  jk+ \lceil{s+k\over 2} \rceil+1 \leq   i\leq s,\,\, i\in \mathcal{Z}\}.
 $$
 \end{definition}

\begin{definition}
Let $$P'(s,k)=L'(s,k)   \cup  R'(s,k) $$
 with the following poset structure:
 \begin{itemize}
  \item If $(a, b), (a', b')\in L'(s,k) $, then $(a, b)$ covers $(a', b')$ if and only if  $|a-a'|\leq {k\over 2}$ and $b=b'+1$;
  \item If $(a, b), (a', b')\in R'(s,k) $, then $(a, b)$ covers $(a', b')$ if and only if  $|a-a'|\leq {k\over 2}$ and $b=b'-1$.
    \end{itemize}
 \end{definition}
For example, the Hasse diagrams of $P'(20,3)$ and $P'(20,4)$ are shown in Figure \ref{P'{20,3}} and Figure \ref{P'{20,4}},   respectively.

 Now we present a map $\chi: P(s, k)\rightarrow P'(s, k)$ as follows.
\begin{itemize}
 \item  If $x=2i-1+2sj$ with $jk+1\leq i\leq \lfloor{s\over 2}\rfloor$ and $j\geq 0$, set  $\chi(x)=(i-{k\over 2}j, j)$;
 \item  If $x=2i-1+2sj$ with  $ jk+ \lceil{s+k\over 2} \rceil+1 \leq   i\leq s$ and $j\geq 0$,  set  $\chi(x)=(s+1-i+{k\over 2}(j+1),-j-1)$.
\end{itemize}
 It is easily seen that the map $\chi: P(s, k)\rightarrow P'(s, k)$ is an isomorphism of posets.
 Furthermore,  $\chi$  maps a nice order ideal of $P(s,k)$ to an order  ideal of $P'(s,k)$ containing    no elements $(a, 0)$ and $(a',-1)$ such that $|a-a'|\leq {k\over 2}$.

 \begin{thm}\label{chi}
 For  given positive integers $s$ and $k$,  the nice order ideals of $P(s, k)$ are in one-to-one correspondence with order  ideals of $P'(s,k)$ containing    no elements $(a, 0)$ and $(a',-1)$ such that $|a-a'|\leq {k\over 2}$
 \end{thm}
 Denote by $J(P'(s,k))$    the set of order  ideals of $P'(s,k)$  containing    no elements $(a, 0)$ and $(a',-1)$ such that $|a-a'|\leq {k\over 2}$.
\subsection{On self-conjugate $(s, s+1, \ldots, s+2k)$-core partitions}

It is routine to check that  the poset $P'(2m,2k)$ is identical with the poset $P'(2m+1, 2k)$.  Relying on Theorems \ref{main1} and \ref{chi}, we are led to the following result.

\begin{thm}\label{mainth2}
For positive integers $m$ and $k$,  self-conjugate $(2m, 2m+1, \ldots, 2m+k)$-core partitions are equinumerous with  self-conjugate $(2m+1, 2m+2, \ldots, 2m+k+1)$-core partitions.
\end{thm}

Combining  Theorems \ref{main1}  and  \ref{mainth2},   we will focus on    the enumeration of   order ideals $I\in J(P'(2m,2k))$  for any positive integers $m$ and $k$. We begin with some definitions and notations.

      If a ballot $(s,k)$-path   is allowed to go below the $x$-axis, then it is called a {\em free} ballot $(s,k)$-path.  Denote by $\mathcal{FB}_n(s,k)$ the set of all free  ballot $(s,k)$-paths of height $n$.
  Let  $\mathcal{SD}(s,k)$ denote the set of all symmetric $(s, k)$-Dyck paths.

Throughout this paper we identify a  (free)  ballot $(s,k)$-path with a word of $U_{k\over 2}$'s, $D_{k\over 2}$'s and $H_{\ell}$'s .    Let $p$ be a step  running from the point $(x_1, y_1)$ to the point $(x_2, y_2)$. Then we say that the point $(x_1,y_1)$ is its  {\em starting
   } point, and the  point $(x_2, y_2)$ is its {\em ending} point. The {\em starting} point and the {\em ending } point  of a path are defined analogously.

 Let $\mathcal{FB}'_n(m,k)$ be the set of free ballot $(m, k)$-paths of height  $n$ which start  with either a horizontal step or a down step.
Let $$\mathcal{Q}(m,k)=(\bigcup_{i=0}^{k-2}\mathcal{FB}'_0(m-i,2k))\cup \mathcal{FB}'_{-1}(m+1,2k).$$

 We describe a map $\varphi:J(P'(2m,2k))\rightarrow \mathcal{Q}(m,k)$ as follows.  Let $I\in J(P'(2m,2k)) $. First we color the  point $(a,b)$ with $a, b\in \mathcal{Z}$ and $0\leq a\leq m+1$ in the plane by black and white according to the following rules.
\begin{itemize}
\item If $b\geq 0$, then
    the point $(a,b)$  is colored by  black if and only   if $(a,b) \in I  $;
    \item If $b<0$, then the point $(a,b)$  is colored by  black if and only    if $(a,b)\notin I  $.
 \end{itemize}
 Now we recursively  define a sequence   $P_0,P_1, \ldots P_n$ of points where  $P_0=(0,0)$. Assume that $  P_{i-1}$ has  been already determined. Let   $P_{i-1}=(a_{i-1}, b_{i-1})$.   We proceed to demonstrate how to get $P_i$.  There are two cases.\\
 Case 1.  There exists a black point   $(a_{i-1}+\ell, b_{i-1})$  with  $1\leq \ell\leq  2k-1$.  If  the point $ (a_{i-1}+k, b_{i-1}+1)$ is black, then  set $P_i=(a_{i-1}+k, b_{i-1}+1)$. Otherwise, we choose $\ell$ to the smallest such integer and set $P_i=(a_{i-1}+\ell, b_{i-1})$.   \\
   Case 2. There does not exist  a black point   $(a_{i-1}+\ell, b_{i-1})$  with  $1\leq \ell\leq  2k-1$. If the point $ (a_{i-1}+k, b_{i-1}-1)$ is black, set $P_i=(a_{i-1}+k, b_{i-1}-1)$. Otherwise, set $P_{i}=P_{i-1}$ and  $n=i-1$.\\
Let $\varphi(I)$ be the resulting  path with lattice points $P_0, P_1, \ldots, P_n$.
For example, let $I=\{(1,0), (2,0), (3,0), (4,0), (5,0), (3,1), (8,-1), (10,-1)\}\in J(P'(20,4))$. By applying the map $\varphi$, we get a  path $\varphi(I)$ as shown in Figure \ref{phi}.

\begin{figure}[h!]
\begin{center}
\begin{picture}(280,160)
\setlength{\unitlength}{8mm}

 \put(0,0){\circle*{0.2}} \put(1,0){\circle*{0.2}} \put(2,0){\circle*{0.2}} \put(3,0){\circle*{0.2}} \put(4,0){\circle*{0.2}} \put(5,0){\circle*{0.2}}
  \put(6,0){\circle*{0.2}} \put(7,0){\circle*{0.2}} \put(8,0){\circle*{0.2}}  \put(9,0){\circle*{0.2}} \put(10,0){\circle*{0.2}} \put(11,0){\circle*{0.2}}

  \put(0,1){\circle*{0.2}} \put(1,1){\circle*{0.2}} \put(2,1){\circle*{0.2}} \put(3,1){\circle*{0.2}} \put(4,1){\circle*{0.2}} \put(5,1){\circle*{0.2}}
  \put(6,1){\circle*{0.2}} \put(7,1){\circle*{0.2}} \put(8,1){\circle*{0.2}}  \put(9,1){\circle*{0.2}} \put(10,1){\circle*{0.2}} \put(11,1){\circle*{0.2}}

\put(0,2){\circle*{0.2}} \put(1,2){\circle*{0.2}} \put(2,2){\circle*{0.2}} \put(3,2){\circle*{0.2}} \put(4,2){\circle*{0.2}} \put(5,2){\circle*{0.2}}
  \put(6,2){\circle*{0.2}} \put(7,2){\circle*{0.2}} \put(8,2){\circle*{0.2}}  \put(9,2){\circle*{0.2}} \put(10,2){\circle*{0.2}} \put(11,2){\circle*{0.2}}

  \put(0,3){\circle*{0.2}} \put(1,3){\circle*{0.2}} \put(2,3){\circle*{0.2}} \put(3,3){\circle*{0.2}} \put(4,3){\circle*{0.2}} \put(5,3){\circle*{0.2}}
  \put(6,3){\circle*{0.2}} \put(7,3){\circle*{0.2}} \put(8,3){\circle{0.2}}  \put(9,3){\circle*{0.2}} \put(10,3){\circle{0.2}} \put(11,3){\circle*{0.2}}

  \put(0,4){\circle{0.2}} \put(1,4){\circle*{0.2}} \put(2,4){\circle*{0.2}} \put(3,4){\circle*{0.2}} \put(4,4){\circle*{0.2}} \put(5,4){\circle*{0.2}}
  \put(6,4){\circle{0.2}} \put(7,4){\circle{0.2}} \put(8,4){\circle{0.2}}  \put(9,4){\circle{0.2}} \put(10,4){\circle{0.2}} \put(11,4){\circle{0.2}}

  \put(0,5){\circle{0.2}} \put(1,5){\circle{0.2}} \put(2,5){\circle{0.2}} \put(3,5){\circle*{0.2}} \put(4,5){\circle{0.2}} \put(5,5){\circle{0.2}}
  \put(6,5){\circle{0.2}} \put(7,5){\circle{0.2}} \put(8,5){\circle{0.2}}  \put(9,5){\circle{0.2}} \put(10,5){\circle{0.2}} \put(11,5){\circle{0.2}}

  \put(0,6){\circle{0.2}} \put(1,6){\circle{0.2}} \put(2,6){\circle{0.2}} \put(3,6){\circle{0.2}} \put(4,6){\circle{0.2}} \put(5,6){\circle{0.2}}
  \put(6,6){\circle{0.2}} \put(7,6){\circle{0.2}} \put(8,6){\circle{0.2}}  \put(9,6){\circle{0.2}} \put(10,6){\circle{0.2}} \put(11,6){\circle{0.2}}

\put(0,4){\line(1,0){1}}\put(1,4){\line(2,1){2}}\put(3,5){\line(2,-1){2}}\put(5,4){\line(2,-1){2}}\put(7,3){\line(2,0){2}} \put(9,3){\line(2,0){2}}

\put(-1.2,3.8){\small$(0,0)$}

\end{picture}
\end{center}
\caption{ The Corresponding path  $\varphi(I)$.}\label{phi}
\end{figure}

\begin{figure}[h!]
\begin{center}
\begin{picture}(280,160)
\setlength{\unitlength}{5mm}

  \put(0,1){\circle{0.2}} \put(1,1){\circle {0.2}} \put(2,1){\circle{0.2}}\put(3,1){\circle{0.2}}\put(4,1){\circle{0.2}}\put(5,1){\circle{0.2}} \put(6,1){\circle{0.2}}\put(7,1){\circle{0.2}}\put(8,1){\circle{0.2}}\put(9,1){\circle{0.2}}\put(10,1){\circle{0.2}}\put(11,1){\circle{0.2}}
  \put(12,1){\circle{0.2}}\put(13,1){\circle{0.2}}\put(14,1){\circle*{0.2}}\put(15,1){\circle*{0.2}}\put(16,1){\circle{0.2}}\put(17,1){\circle{0.2}}
  \put(18,1){\circle{0.2}}\put(19,1){\circle{0.2}}\put(20,1){\circle{0.2}}\put(21,1){\circle{0.2}}\put(22,1){\circle*{0.2}}  \put(23,1){\circle{0.2}}
   \put(24,1){\circle{0.2}}

  \put(0,0){\circle{0.2}} \put(1,0){\circle {0.2}} \put(2,0){\circle{0.2}}\put(3,0){\circle{0.2}}\put(4,0){\circle{0.2}}\put(5,0){\circle{0.2}} \put(6,0){\circle{0.2}}\put(7,0){\circle{0.2}}\put(8,0){\circle{0.2}}\put(9,0){\circle{0.2}}\put(10,0){\circle{0.2}}\put(11,0){\circle{0.2}}
  \put(12,0){\circle{0.2}}\put(13,0){\circle{0.2}}\put(14,0){\circle{0.2}}\put(15,0){\circle{0.2}}\put(16,0){\circle{0.2}}\put(17,0){\circle*{0.2}}
  \put(18,0){\circle{0.2}}\put(19,0){\circle*{0.2}}\put(20,0){\circle*{0.2}}\put(21,0){\circle{0.2}}\put(22,0){\circle{0.2}} \put(23,0){\circle{0.2}}
   \put(24,0){\circle{0.2}}

  \put(0,2){\circle*{0.2}} \put(1,2){\circle {0.2}} \put(2,2){\circle{0.2}}\put(3,2){\circle*{0.2}}\put(4,2){\circle{0.2}}\put(5,2){\circle{0.2}} \put(6,2){\circle{0.2}}\put(7,2){\circle{0.2}}\put(8,2){\circle{0.2}}\put(9,2){\circle{0.2}}\put(10,2){\circle{0.2}}\put(11,2){\circle{0.2}}
  \put(12,2){\circle*{0.2}}\put(13,2){\circle{0.2}}\put(14,2){\circle{0.2}}\put(15,2){\circle{0.2}}\put(16,2){\circle{0.2}}\put(17,2){\circle{0.2}}
  \put(18,2){\circle{0.2}}\put(19,2){\circle{0.2}}\put(20,2){\circle{0.2}}\put(21,2){\circle{0.2}}\put(22,2){\circle{0.2}} \put(23,2){\circle{0.2}}
   \put(24,2){\circle*{0.2}}

  \put(0,3){\circle{0.2}} \put(1,3){\circle {0.2}} \put(2,3){\circle{0.2}}\put(3,3){\circle{0.2}}\put(4,3){\circle{0.2}}\put(5,3){\circle*{0.2}} \put(6,3){\circle*{0.2}}\put(7,3){\circle{0.2}}\put(8,3){\circle{0.2}}\put(9,3){\circle{0.2}}\put(10,3){\circle*{0.2}}\put(11,3){\circle{0.2}}
  \put(12,3){\circle{0.2}}\put(13,3){\circle{0.2}}\put(14,3){\circle{0.2}}\put(15,3){\circle{0.2}}\put(16,3){\circle{0.2}}\put(17,3){\circle{0.2}}
  \put(18,3){\circle{0.2}}\put(19,3){\circle{0.2}}\put(20,3){\circle{0.2}}\put(21,3){\circle{0.2}}\put(22,3){\circle{0.2}} \put(23,3){\circle{0.2}}
   \put(24,3){\circle{0.2}}
   \put(0,4){\circle{0.2}} \put(1,4){\circle {0.2}} \put(2,4){\circle{0.2}}\put(3,4){\circle{0.2}}\put(4,4){\circle{0.2}}\put(5,4){\circle{0.2}} \put(6,4){\circle{0.2}}\put(7,4){\circle{0.2}}\put(8,4){\circle*{0.2}}\put(9,4){\circle{0.2}}\put(10,4){\circle{0.2}}\put(11,4){\circle{0.2}}
  \put(12,4){\circle{0.2}}\put(13,4){\circle{0.2}}\put(14,4){\circle{0.2}}\put(15,4){\circle{0.2}}\put(16,4){\circle{0.2}}\put(17,4){\circle{0.2}}
  \put(18,4){\circle{0.2}}\put(19,4){\circle{0.2}}\put(20,4){\circle{0.2}}\put(21,4){\circle{0.2}}\put(22,4){\circle{0.2}} \put(23,4){\circle{0.2}}
   \put(24,4){\circle{0.2}}
\put(0,2){\line(3,0){3}} \put(3,2){\line(2,1){2}}  \put(5,3){\line(1,0){1}} \put(6,3){\line(2,1){2}} \put(8,4){\line(2,-1){2}} \put(10,3){\line(2,-1){2}}
\put(12,2){\line(2,-1){2}} \put(14,1){\line(1,0){1}} \put(15,1){\line(2,-1){2}} \put(17,0){\line(2,0){2}}  \put(19,0){\line(1,0){1}}  \put(20,0){\line(2,1){2}}
\put(22,1){\line(2,1){2}}

\put(-2,1.8){\small$(0,0)$} \put(17,-0.5){\small$a$}

\end{picture}
\end{center}
\caption{ The decomposition  of a path $ P=H_3U_2H_1U_2D_2D_2D_2H_1D_2H_2H_1U_2U_2$.}\label{alpha1}
\end{figure}

\begin{figure}[h!]
\begin{center}
\begin{picture}(280,160)
\setlength{\unitlength}{5mm}

  \put(0,1){\circle{0.2}} \put(1,1){\circle {0.2}} \put(2,1){\circle{0.2}}\put(3,1){\circle{0.2}}\put(4,1){\circle{0.2}}\put(5,1){\circle*{0.2}} \put(6,1){\circle{0.2}}\put(7,1){\circle{0.2}}\put(8,1){\circle{0.2}}\put(9,1){\circle{0.2}}\put(10,1){\circle{0.2}}\put(11,1){\circle{0.2}}
  \put(12,1){\circle{0.2}}\put(13,1){\circle{0.2}}\put(14,1){\circle{0.2}}\put(15,1){\circle{0.2}}\put(16,1){\circle{0.2}}\put(17,1){\circle{0.2}}
  \put(18,1){\circle{0.2}}\put(19,1){\circle{0.2}}\put(20,1){\circle{0.2}}\put(21,1){\circle{0.2}}\put(22,1){\circle{0.2}}  \put(23,1){\circle{0.2}}
   \put(24,1){\circle{0.2}}

   \put(0,0){\circle*{0.2}} \put(1,0){\circle {0.2}} \put(2,0){\circle*{0.2}}\put(3,0){\circle*{0.2}}\put(4,0){\circle{0.2}}\put(5,0){\circle{0.2}} \put(6,0){\circle{0.2}}\put(7,0){\circle{0.2}}\put(8,0){\circle{0.2}}\put(9,0){\circle{0.2}}\put(10,0){\circle{0.2}}\put(11,0){\circle{0.2}}
  \put(12,0){\circle{0.2}}\put(13,0){\circle{0.2}}\put(14,0){\circle{0.2}}\put(15,0){\circle{0.2}}\put(16,0){\circle{0.2}}\put(17,0){\circle{0.2}}
  \put(18,0){\circle{0.2}}\put(19,0){\circle{0.2}}\put(20,0){\circle{0.2}}\put(21,0){\circle{0.2}}\put(22,0){\circle{0.2}}  \put(23,0){\circle{0.2}}
   \put(24,0){\circle{0.2}}

    \put(0,2){\circle{0.2}} \put(1,2){\circle {0.2}} \put(2,2){\circle{0.2}}\put(3,2){\circle{0.2}}\put(4,2){\circle{0.2}}\put(5,2){\circle{0.2}} \put(6,2){\circle{0.2}}\put(7,2){\circle*{0.2}}\put(8,2){\circle{0.2}}\put(9,2){\circle{0.2}}\put(10,2){\circle{0.2}}\put(11,2){\circle{0.2}}
  \put(12,2){\circle{0.2}}\put(13,2){\circle{0.2}}\put(14,2){\circle{0.2}}\put(15,2){\circle{0.2}}\put(16,2){\circle{0.2}}\put(17,2){\circle{0.2}}
  \put(18,2){\circle{0.2}}\put(19,2){\circle{0.2}}\put(20,2){\circle{0.2}}\put(21,2){\circle{0.2}}\put(22,2){\circle{0.2}}  \put(23,2){\circle{0.2}}
   \put(24,2){\circle{0.2}}

    \put(0,3){\circle{0.2}} \put(1,3){\circle {0.2}} \put(2,3){\circle{0.2}}\put(3,3){\circle{0.2}}\put(4,3){\circle{0.2}}\put(5,3){\circle{0.2}} \put(6,3){\circle{0.2}}\put(7,3){\circle{0.2}}\put(8,3){\circle{0.2}}\put(9,3){\circle*{0.2}}\put(10,3){\circle*{0.2}}\put(11,3){\circle{0.2}}
  \put(12,3){\circle{0.2}}\put(13,3){\circle{0.2}}\put(14,3){\circle{0.2}}\put(15,3){\circle{0.2}}\put(16,3){\circle{0.2}}\put(17,3){\circle{0.2}}
  \put(18,3){\circle{0.2}}\put(19,3){\circle{0.2}}\put(20,3){\circle{0.2}}\put(21,3){\circle{0.2}}\put(22,3){\circle{0.2}}  \put(23,3){\circle{0.2}}
   \put(24,3){\circle{0.2}}

   \put(0,4){\circle{0.2}} \put(1,4){\circle {0.2}} \put(2,4){\circle{0.2}}\put(3,4){\circle{0.2}}\put(4,4){\circle{0.2}}\put(5,4){\circle{0.2}} \put(6,4){\circle{0.2}}\put(7,4){\circle{0.2}}\put(8,4){\circle{0.2}}\put(9,4){\circle{0.2}}\put(10,4){\circle{0.2}}\put(11,4){\circle{0.2}}
  \put(12,4){\circle*{0.2}}\put(13,4){\circle{0.2}}\put(14,4){\circle{0.2}}\put(15,4){\circle{0.2}}\put(16,4){\circle{0.2}}\put(17,4){\circle{0.2}}
  \put(18,4){\circle{0.2}}\put(19,4){\circle{0.2}}\put(20,4){\circle{0.2}}\put(21,4){\circle*{0.2}}\put(22,4){\circle{0.2}}  \put(23,4){\circle{0.2}}
   \put(24,4){\circle*{0.2}}

   \put(0,5){\circle{0.2}} \put(1,5){\circle {0.2}} \put(2,5){\circle{0.2}}\put(3,5){\circle{0.2}}\put(4,5){\circle{0.2}}\put(5,5){\circle{0.2}} \put(6,5){\circle{0.2}}\put(7,5){\circle{0.2}}\put(8,5){\circle{0.2}}\put(9,5){\circle{0.2}}\put(10,5){\circle{0.2}}\put(11,5){\circle{0.2}}
  \put(12,5){\circle{0.2}}\put(13,5){\circle{0.2}}\put(14,5){\circle*{0.2}}\put(15,5){\circle{0.2}}\put(16,5){\circle{0.2}}\put(17,5){\circle{0.2}}
  \put(18,5){\circle*{0.2}}\put(19,5){\circle*{0.2}}\put(20,5){\circle{0.2}}\put(21,5){\circle{0.2}}\put(22,5){\circle{0.2}}  \put(23,5){\circle{0.2}}
   \put(24,5){\circle{0.2}}

   \put(0,6){\circle{0.2}} \put(1,6){\circle {0.2}} \put(2,6){\circle{0.2}}\put(3,6){\circle{0.2}}\put(4,6){\circle{0.2}}\put(5,6){\circle{0.2}} \put(6,6){\circle{0.2}}\put(7,6){\circle{0.2}}\put(8,6){\circle{0.2}}\put(9,6){\circle{0.2}}\put(10,6){\circle{0.2}}\put(11,6){\circle{0.2}}
  \put(12,6){\circle{0.2}}\put(13,6){\circle{0.2}}\put(14,6){\circle{0.2}}\put(15,6){\circle{0.2}}\put(16,6){\circle*{0.2}}\put(17,6){\circle{0.2}}
  \put(18,6){\circle{0.2}}\put(19,6){\circle{0.2}}\put(20,6){\circle{0.2}}\put(21,6){\circle{0.2}}\put(22,6){\circle{0.2}}  \put(23,6){\circle{0.2}}
   \put(24,6){\circle{0.2}}

\put(0,0){\line(2,0){2}} \put(2,0){\line(1,0){1}} \put(3,0){\line(2,1){2}} \put(5,1){\line(2,1){2}} \put(7,2){\line(2,1){2}}\put(9,3){\line(1,0){1}}
\put(10,3){\line(2,1){2}}\put(12,4){\line(2,1){2}}\put(14,5){\line(2,1){2}} \put(16,6){\line(2,-1){2}} \put(18,5){\line(1,0){1}}
\put(19,5){\line(2,-1){2}} \put(21,4){\line(3,0){3}}

\put(-2,-0.2){\small$(0,0)$}  \put(7.2,1.5){\small$b$}

\end{picture}
\end{center}
\caption{ The decomposition of a path $ Q=H_2H_1U_2U_2U_2H_1U_2U_2U_2D_2H_1D_2H_3$.}\label{alpha2}
\end{figure}
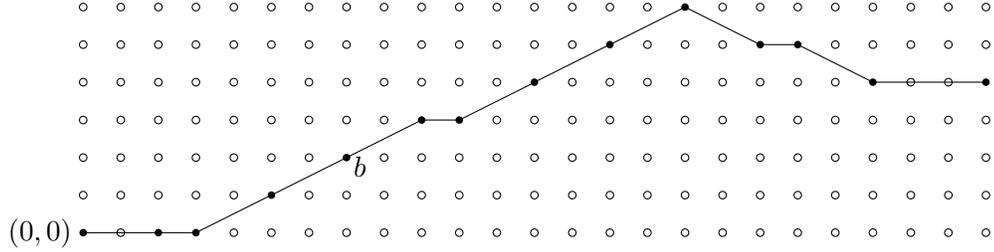

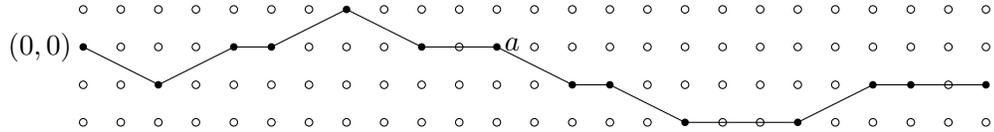
\begin{figure}[h!]
\begin{center}
\begin{picture}(280,160)
\setlength{\unitlength}{5mm}

  \put(0,1){\circle{0.2}} \put(1,1){\circle {0.2}} \put(2,1){\circle*{0.2}}\put(3,1){\circle{0.2}}\put(4,1){\circle{0.2}}\put(5,1){\circle{0.2}} \put(6,1){\circle{0.2}}\put(7,1){\circle{0.2}}\put(8,1){\circle{0.2}}\put(9,1){\circle{0.2}}\put(10,1){\circle{0.2}}\put(11,1){\circle{0.2}}
  \put(12,1){\circle{0.2}}\put(13,1){\circle*{0.2}}\put(14,1){\circle*{0.2}}\put(15,1){\circle{0.2}}\put(16,1){\circle{0.2}}\put(17,1){\circle{0.2}}
  \put(18,1){\circle{0.2}}\put(19,1){\circle{0.2}}\put(20,1){\circle{0.2}}\put(21,1){\circle*{0.2}}\put(22,1){\circle*{0.2}}  \put(23,1){\circle{0.2}}
   \put(24,1){\circle*{0.2}}

  \put(0,0){\circle{0.2}} \put(1,0){\circle {0.2}} \put(2,0){\circle{0.2}}\put(3,0){\circle{0.2}}\put(4,0){\circle{0.2}}\put(5,0){\circle{0.2}} \put(6,0){\circle{0.2}}\put(7,0){\circle{0.2}}\put(8,0){\circle{0.2}}\put(9,0){\circle{0.2}}\put(10,0){\circle{0.2}}\put(11,0){\circle{0.2}}
  \put(12,0){\circle{0.2}}\put(13,0){\circle{0.2}}\put(14,0){\circle{0.2}}\put(15,0){\circle{0.2}}\put(16,0){\circle*{0.2}}\put(17,0){\circle{0.2}}
  \put(18,0){\circle{0.2}}\put(19,0){\circle*{0.2}}\put(20,0){\circle{0.2}}\put(21,0){\circle{0.2}}\put(22,0){\circle{0.2}}  \put(23,0){\circle{0.2}}
   \put(24,0){\circle{0.2}}

   \put(0,2){\circle*{0.2}} \put(1,2){\circle {0.2}} \put(2,2){\circle{0.2}}\put(3,2){\circle{0.2}}\put(4,2){\circle*{0.2}}\put(5,2){\circle*{0.2}} \put(6,2){\circle{0.2}}\put(7,2){\circle{0.2}}\put(8,2){\circle{0.2}}\put(9,2){\circle*{0.2}}\put(10,2){\circle{0.2}}\put(11,2){\circle*{0.2}}
  \put(12,2){\circle{0.2}}\put(13,2){\circle{0.2}}\put(14,2){\circle{0.2}}\put(15,2){\circle{0.2}}\put(16,2){\circle{0.2}}\put(17,2){\circle{0.2}}
  \put(18,2){\circle{0.2}}\put(19,2){\circle{0.2}}\put(20,2){\circle{0.2}}\put(21,2){\circle{0.2}}\put(22,2){\circle{0.2}}  \put(23,2){\circle{0.2}}
   \put(24,2){\circle{0.2}}

    \put(0,3){\circle{0.2}} \put(1,3){\circle {0.2}} \put(2,3){\circle{0.2}}\put(3,3){\circle{0.2}}\put(4,3){\circle{0.2}}\put(5,3){\circle{0.2}} \put(6,3){\circle{0.2}}\put(7,3){\circle*{0.2}}\put(8,3){\circle{0.2}}\put(9,3){\circle{0.2}}\put(10,3){\circle{0.2}}\put(11,3){\circle{0.2}}
  \put(12,3){\circle{0.2}}\put(13,3){\circle{0.2}}\put(14,3){\circle{0.2}}\put(15,3){\circle{0.2}}\put(16,3){\circle{0.2}}\put(17,3){\circle{0.2}}
  \put(18,3){\circle{0.2}}\put(19,3){\circle{0.2}}\put(20,3){\circle{0.2}}\put(21,3){\circle{0.2}}\put(22,3){\circle{0.2}}  \put(23,3){\circle{0.2}}
   \put(24,3){\circle{0.2}}

\put(0,2){\line(2,-1){2}}  \put(2,1){\line(2,1){2}}  \put(4,2){\line(1,0){1}}  \put(5,2){\line(2,1){2}}  \put(7,3){\line(2,-1){2}}   \put(9,2){\line(2,0){2}}
\put(11,2){\line(2,-1){2}}  \put(13,1){\line(1,0){1}}   \put(14,1){\line(2,-1){2}}  \put(16,0){\line(3,0){3}}  \put(19,0){\line(2,1){2}} \put(21,1){\line(1,0){1}} \put(22,1){\line(2,0){2}}

\put(-2,1.8){\small$(0,0)$}  \put(11.2,1.9){\small$a$}

\end{picture}
\end{center}
\caption{ The decomposition of a path $P=D_2U_2H_1U_2D_2H_2D_2H_1D_2H_3U_2H_1H_2$.}\label{beta1}
\end{figure}

\begin{figure}[h!]
\begin{center}
\begin{picture}(280,160)
\setlength{\unitlength}{5mm}

  \put(0,1){\circle{0.2}} \put(1,1){\circle {0.2}} \put(2,1){\circle*{0.2}}\put(3,1){\circle*{0.2}}\put(4,1){\circle{0.2}}\put(5,1){\circle{0.2}} \put(6,1){\circle{0.2}}\put(7,1){\circle*{0.2}}\put(8,1){\circle{0.2}}\put(9,1){\circle*{0.2}}\put(10,1){\circle{0.2}}\put(11,1){\circle{0.2}}
  \put(12,1){\circle{0.2}}\put(13,1){\circle {0.2}}\put(14,1){\circle {0.2}}\put(15,1){\circle{0.2}}\put(16,1){\circle{0.2}}\put(17,1){\circle{0.2}}
  \put(18,1){\circle{0.2}}\put(19,1){\circle{0.2}}\put(20,1){\circle{0.2}}\put(21,1){\circle{0.2}}\put(22,1){\circle{0.2}}  \put(23,1){\circle{0.2}}
   \put(24,1){\circle{0.2}}

  \put(0,0){\circle*{0.2}} \put(1,0){\circle {0.2}} \put(2,0){\circle{0.2}}\put(3,0){\circle{0.2}}\put(4,0){\circle{0.2}}\put(5,0){\circle{0.2}} \put(6,0){\circle{0.2}}\put(7,0){\circle{0.2}}\put(8,0){\circle{0.2}}\put(9,0){\circle{0.2}}\put(10,0){\circle{0.2}}\put(11,0){\circle{0.2}}
  \put(12,0){\circle{0.2}}\put(13,0){\circle{0.2}}\put(14,0){\circle{0.2}}\put(15,0){\circle{0.2}}\put(16,0){\circle{0.2}}\put(17,0){\circle{0.2}}
  \put(18,0){\circle{0.2}}\put(19,0){\circle{0.2}}\put(20,0){\circle{0.2}}\put(21,0){\circle{0.2}}\put(22,0){\circle{0.2}}  \put(23,0){\circle{0.2}}
   \put(24,0){\circle{0.2}}

   \put(0,2){\circle{0.2}} \put(1,2){\circle {0.2}} \put(2,2){\circle{0.2}}\put(3,2){\circle{0.2}}\put(4,2){\circle{0.2}}\put(5,2){\circle*{0.2}} \put(6,2){\circle{0.2}}\put(7,2){\circle{0.2}}\put(8,2){\circle{0.2}}\put(9,2){\circle{0.2}}\put(10,2){\circle{0.2}}\put(11,2){\circle*{0.2}}
  \put(12,2){\circle{0.2}}\put(13,2){\circle{0.2}}\put(14,2){\circle{0.2}}\put(15,2){\circle{0.2}}\put(16,2){\circle{0.2}}\put(17,2){\circle{0.2}}
  \put(18,2){\circle{0.2}}\put(19,2){\circle{0.2}}\put(20,2){\circle{0.2}}\put(21,2){\circle{0.2}}\put(22,2){\circle{0.2}}  \put(23,2){\circle{0.2}}
   \put(24,2){\circle{0.2}}

    \put(0,3){\circle{0.2}} \put(1,3){\circle {0.2}} \put(2,3){\circle{0.2}}\put(3,3){\circle{0.2}}\put(4,3){\circle{0.2}}\put(5,3){\circle{0.2}} \put(6,3){\circle{0.2}}\put(7,3){\circle{0.2}}\put(8,3){\circle{0.2}}\put(9,3){\circle{0.2}}\put(10,3){\circle{0.2}}\put(11,3){\circle{0.2}}
  \put(12,3){\circle{0.2}}\put(13,3){\circle*{0.2}}\put(14,3){\circle*{0.2}}\put(15,3){\circle{0.2}}\put(16,3){\circle{0.2}}\put(17,3){\circle{0.2}}
  \put(18,3){\circle{0.2}}\put(19,3){\circle{0.2}}\put(20,3){\circle{0.2}}\put(21,3){\circle*{0.2}}\put(22,3){\circle*{0.2}}  \put(23,3){\circle{0.2}}
   \put(24,3){\circle*{0.2}}

\put(0,4){\circle{0.2}} \put(1,4){\circle {0.2}} \put(2,4){\circle{0.2}}\put(3,4){\circle{0.2}}\put(4,4){\circle{0.2}}\put(5,4){\circle{0.2}} \put(6,4){\circle{0.2}}\put(7,4){\circle{0.2}}\put(8,4){\circle{0.2}}\put(9,4){\circle{0.2}}\put(10,4){\circle{0.2}}\put(11,4){\circle{0.2}}
  \put(12,4){\circle{0.2}}\put(13,4){\circle{0.2}}\put(14,4){\circle{0.2}}\put(15,4){\circle{0.2}}\put(16,4){\circle*{0.2}}\put(17,4){\circle{0.2}}
  \put(18,4){\circle{0.2}}\put(19,4){\circle*{0.2}}\put(20,4){\circle{0.2}}\put(21,4){\circle{0.2}}\put(22,4){\circle{0.2}}  \put(23,4){\circle{0.2}}
   \put(24,4){\circle{0.2}}

\put(0,0){\line(2,1){2}}   \put(2,1){\line(1,0){1}}   \put(3,1){\line(2,1){2}}   \put(5,2){\line(2,-1){2}}  \put(7,1){\line(2,0){2}}  \put(9,1){\line(2,1){2}}
\put(11,2){\line(2,1){2}}   \put(13,3){\line(1,0){1}}  \put(14,3){\line(2,1){2}}   \put(16,4){\line(3,0){3}}   \put(19,4){\line(2,-1){2}}  \put(21,3){\line(1,0){1}} \put(22,3){\line(2,0){2}}

\put(-2,-0.2){\small$(0,0)$}   \put(11.2,1.6){\small$b$}

\end{picture}
\end{center}
\caption{ The decomposition of a path $Q=U_2H_1U_2D_2H_2U_2U_2H_1U_2H_3D_2H_1H_2$.}\label{beta2}
\end{figure}
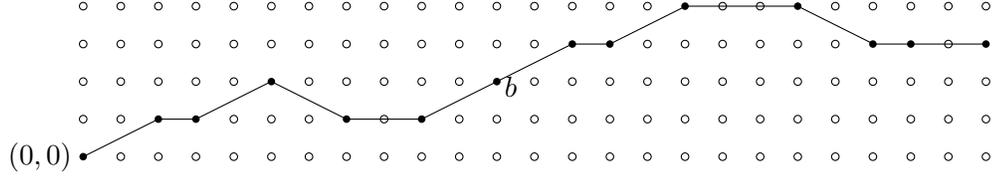

\begin{lem}\label{lemphi3}
For any order ideal $I\in J(P'(2m,2k))$,  the resulting path $\varphi(I)$ verifies the following properties.
\begin{itemize}
 \item[{\upshape (\rmnum{1})}] The points   lying     above the path $\varphi(I)$ are white;
\item[ {\upshape (\rmnum{2})}] The points   lying     below the path $\varphi(I)$ are black.
\end{itemize}
 \end{lem}
\pf Here we only verify  {\upshape (\rmnum{1})}.  By similar arguments,  one can prove {\upshape (\rmnum{2})}. We follow notations in the definition of $\varphi$.     Let $P^{(\ell)}$ be  the partial path of $\varphi(I)$  whose lattice points are $P_0, P_1, \ldots, P_{\ell}$.  We prove   by induction on $\ell$.   Clearly the assertion holds for $\ell=0$.  Assume that  the points   lying    above the path $P^{(\ell-1)} $ are white.  Let $P_{\ell-1}=(a,b)$.   Now we proceed to show that the assertion also holds for $P^{(\ell)}$.
 We have three cases.

 Case 1. $P_{\ell}=(a+k, b+1)$.     When $y>b+1$ and $a<x\leq a+k$,  the point $(x,y)$ must be white. If not, by the cover relation of $I$, the point $(a, y-1)$  is also black, which contradicts the induction hypothesis. It remains to show that the point $(x,b+1)$ is white for any $a<x<a+k$. If not,  Let $c$ be the smallest integer such that  the point $(c,b+1)$ is black and $a<c<a+k$.     By  the cover relation of $I$, it follows that  the point $(c-k, b)$ is black.  By the induction hypothesis, we must have $P_{\ell-2}=(a-k, b+1)$.  Then by the construction of the path, we must  have $P_{\ell-1}=(c,b+1)$, which  contradicts the selection of $P_{\ell-1}$.

 Case 2.  $P_{\ell}=(a+i, b)$ for some $1\leq i\leq 2k-1$.  We claim that the point $(x,b+1)$ is white for any $a<x\leq a+i$. If not,  Let $c$ be the smallest integer such that  the point $(c,b+1)$ is black and $a<c\leq a+i$.  According to the selection of $P_\ell$,  we have $c\neq a+k$.  By the cover relation of $I$, the point $(c-k, b)$ is black.  Hence,  by the induction hypothesis, we must have $P_{\ell-2}=(a-k, b+1)$.  Then by the construction of the path, we must  have $P_{\ell-1}=(c,b+1)$, which  contradicts the selection of $P_{\ell-1}$.   Hence,  by the cover relation of $I$, the point $(x,y)$ must be white for any $y>b+1$ and $a<x\leq a+i$.

 Case 3.  $P_{\ell}=(a+k, b-1)$.  Due to the selection of the lattice point $P_{\ell}$,  the point $(x,b)$ is white for any $a<x\leq  a+k$. Then by the cover relation of $I$,    the point $(x,y)$ must be white for any $y>b+1$ and $a<x\leq a+k$.  This completes the proof. \qed

   \begin{lem}\label{lemphi1}
 The map $\varphi$ is well defined, that is, for any $I\in J(P'(2m,2k))$  we have   $ \varphi(I)\in\mathcal{Q}(m,k)$.
\end{lem}
\pf Let $\varphi(I)$ be   path with lattice points $P_0, P_1, \ldots, P_n$.
   It is apparent that the resulting path is a lattice path  consisting of up steps $U_k$, down steps $D_k$ and horizontal steps $H_\ell$ for some $1\leq \ell\leq 2k-1$. Moreover, since $(k, 1)$ is not contained in $P'(2m,2k)$, $\varphi(I)$ must start with either a down step or a horizontal step. In order to show that $\varphi$ is well defined, it remains  to verify that we have either $P_{n}=(m-i, 0)$ for some $0\leq i\leq k-2$  or $P_n=(m+1, -1)$. Notice that $(k,-1)$ is not contained in  $P'(m,k)$. Thus the point $(k,-1)$ is always colored by black. This ensures that we have $n>0$.

 Suppose that $P_n=(a,b)$. According to the construction of the path,  the point $P_n$ is black.  If  $b> 0$,  then we have $(a,b)\in I$ and $(a+k, b-1)\in I$ by the cover relation of $I$.  It follows that  the point $(a+k, b-1)$ is black  and  the path can not terminate at the   point $(a,b)$.  So we have $b\leq 0$.

 We claim that $P_n\neq (m-i, 0)$ for any $i\geq k-1$. If not, we have $(m-i, 0)\in I$, which implies that   $(m-i+k, -1)\notin I$.  It follows  that the point $(m-i+k,-1)$ is black  and the path can not terminate at the   point $(a,b)$.

 Suppose that $b<0$ and $a\leq m+1-k$.     Since the point $(a, b)$ is the ending point of $\varphi(I)$, the point $(a+k, b-1) $ must be  white.  This implies that we have $(a+k,b-1)\in I$. Recall that  $ (a+k,b-1) $ covers $(a,b)$. Therefore,   we have $(a, b)\in I$ and  the point   $(a,b)$ is white,   a contradiction.

 Suppose that $b<0$ and $  m+2-k\leq a\leq m$.   Since $m+1-a\leq k-1$ and $\varphi(I)$ terminates at $(a,b)$, the point   $(m+1,b)$ must be  white.  By the cover relation of $I$, we have $(m+1,-1)\in I$,    a contradiction.

 We claim that $P_n\neq (m+1, b)$ for any $b<-1$. If not,  by {\upshape (\rmnum{1})} of Lemma \ref{lemphi3}, the point $(m+1, -1)$ is white,  yielding a contradiction with the fact that the point $(m+1, -1)$ is black.

    So far, we have reached the conclusion that   we have either $P_{n}=(m-i, 0)$ for some $0\leq i\leq k-2$  or $P_n=(m+1, -1)$. This implies that the map $\varphi$ is well defined, completing the  proof. \qed

\begin{lem}\label{lemphi2}
For any order ideal $I\in J(P'(2m,2k))$,  the resulting path $\varphi(I)$ verifies the following properties.
\begin{itemize}
\item[{\upshape (\rmnum{1})}] The points   on  the path $\varphi(I)$  are black if and only if they are ending points of steps;
\item[ {\upshape (\rmnum{2})}] For $b<0$, the point $(a,b)$       to the right of the ending point of  the path $\varphi(I)$ is black;
\item[ {\upshape (\rmnum{3})}] For $b\geq 0$, the point $(a,b)$        to the right of the ending point of  the path $\varphi(I)$ is white.
 \end{itemize}
\end{lem}

\pf  Point  {\upshape (\rmnum{1})} follows  directly from the  construction of  $\varphi(I)$.

Recall that the ending point of $\varphi(I)$ is either the point $(m-i, 0)$ for some $0\leq i\leq k-2$ or the point $(m+1, -1)$. Clearly, it suffices to prove  $\varphi(I)$ verifies properties    {\upshape (\rmnum{2})} and  {\upshape (\rmnum{3})} for the former case.  Since $\varphi(I)$ terminates at $(m-i,0)$, the point  $(a, 0)$ is white for any $m-i<a\leq m+1$. Hence by the cover relation of $I$, the point $(a,b)$ is white for any $m-i<a\leq m+1$ and $b\geq 0$.  Moreover,  the point  $(a, -1)$ is black for any $m-i<a\leq m+1$. This implies that  the point $(a,b)$ is black for any $m-i<a\leq m+1$ and $b\leq -1$. Thus, $\varphi(I)$ has  properties    {\upshape (\rmnum{2})} and  {\upshape (\rmnum{3})}. This completes the proof. \qed

\begin{thm}\label{mainth3}
The map $\varphi$ is a bijection between  $J(P'(2m,2k)) $ and $\mathcal{Q}(m,k)$.
\end{thm}
\pf   First we give  the inverse map $\varphi'$ from the set $\mathcal{Q}(m,k)$ to the set $J(P'(2m,2k))$. Let $P\in \mathcal{Q}(m,k)$. First we color the  point $(a,b)$ with $a, b\in \mathcal{Z}$ and $0\leq a\leq m+1$ in the plane by black and white according to the following rules.
\begin{itemize}
\item If the point $(a,b)$ lies below the path, we color it by black;
    \item  If the point $(a,b)$ lies above the path, we color it by white;
    \item   For a point $(a,b)$   on the path, we color it by black if and only if it is  an    ending point  of a step;
        \item For  the point (a,b)   to the right to the endpoint of   $P$, we color it by white when $b\geq 0$,  and color it by black when $b<0$.
  \end{itemize}
For $b<0$,  let $(a,b)\in \varphi'(P)$  if and only if the point $(a,b)$ is white. For $b\geq 0$, let $(a,b)\in \varphi'(P)$ if and only if the point $(a,b)$ is black.  It is easily seen  that $\varphi'(P) \in J(P'(2m,2k))$.
Lemmas \ref{lemphi3} and \ref{lemphi2} guarantee  that the maps $\varphi $ and $\varphi'$ are inverses of each other, and thus the map $\varphi$ is a bijection. This completes the proof. \qed

In the following, we aim to establish a bijection between $\mathcal{Q}(m,k)$ and $\mathcal{SD}(2m,2k)$.

Let $\mathcal{B}^e(s,k)$ and $\mathcal{B}^o(s,k)$ denote the set of ballot $(s,k)$-paths of even height and the set of ballot $(s,k)$-paths of odd height, respectively.  The {\em complement } of the path $P$, denoted by $   P^c$,  is  obtained from $P$ by interchanging the up steps and the down steps.   The {\em reverse complement } of the path $P$,
denoted by $ P^{^{rc}}$,  is obtained from $P$ by   encoding the steps of $P$ from right to left and  interchanging the up steps and the down steps.
For example, the complement  of the path  $P=H_1U_3D_3D_3U_3D_3H_2$ is given by  $H_1D_3U_3U_3D_3U_3H_2  $, while the reverse complement of the path $P$ is given by $H_2U_3D_3U_3U_3D_3H_1$.

\begin{thm} \label{alpha}
There is a bijection $\alpha$  between  the set $\mathcal{FB}_0(s,k)$ and the set $\mathcal{B}^{e}(s,k)$.
\end{thm}
\pf  Let $P\in \mathcal{FB}_0(s,k) $. Suppose that $a$ is  the lowest point of $P$. If there are more than one such lowest point, we choose $a$ to be the
leftmost one. Now we proceed to construct a path $\alpha(P)$.
    Clearly, $P$ can be uniquely decomposed as  $P=P^{(1)}P^{(2)}$,
    where $P^{(1)}$ is the subpath of $P$ which goes from the starting point of $P$ to the point $a$, and $P^{(2)}$ is the remaining subpath of $P$.
          Set $\alpha(P)=P^{(2)}   ( P^{(1)})^{rc}$.
 It is apparent that  the resulting path  is a ballot $(s,k)$-path of  even height.

 In order to show that the map $\alpha$ is a bijection, we describe a map $\alpha':\mathcal{B}^{e}(s,k) \rightarrow \mathcal{FB}_0(s,k)$. Let $Q $ be a ballot $(s,k)$-path of height $2\ell$ for some $\ell\geq 0$. Now we proceed to demonstrate how to get a path $\alpha'(Q)$.  Let $b$ the rightmost lattice point of $Q$ which is on the line $y=\ell$. Then $Q$ can be uniquely decomposed as $Q=Q^{(1)}Q^{(2)}$, where $Q^{(1)}$ is the subpath of $Q$ which goes from the starting point of $Q$ to the point $b$, and $Q^{(2)}$ is the remaining subpath of $Q$. Set $\alpha'(Q)= (Q^{(2)})^{rc} Q^{(1)}$. One can easily verify that the resulting path $\alpha'(Q)\in \mathcal{FB}_0(s,k)$. Moreover, the point $b$ becomes the leftmost  lowest point of the resulting path.  This ensures that the maps $\alpha$ and $\alpha'$ are inverses of each other. Thus the map $\alpha$ is a bijection, completing the proof. \qed

\begin{example}
The decomposition of a path $$P=H_3U_2H_1U_2D_2D_2D_2H_1D_2H_2H_1U_2U_2$$ is shown in Figure \ref{alpha1}. By applying the map $\alpha$, its corresponding path is illustrated in Figure \ref{alpha2}.
\end{example}

\begin{example}
The decomposition of a path $$Q=H_2H_1U_2U_2U_2H_1U_2U_2U_2D_2H_1D_2H_3$$ is shown in Figure \ref{alpha2}. By applying the map $\alpha'$, its corresponding path is illustrated in Figure \ref{alpha1}.
\end{example}

\begin{thm}\label{beta}
There is a bijection $\beta$ between  the set $\mathcal{FB}_{-1}(s,k)$ and the set $\mathcal{B}^{o}(s,k)$.
\end{thm}
 \pf  Let $P\in \mathcal{FB}_{-1}(s,k) $.
      Now we proceed to construct a path $\beta(P)$ as follows. Suppose that $a$ is  the rightmost lattice  point of $P$ which is on the $x$-axis.
    Then $P$ can be uniquely decomposed as  $P=P^{(1)}P^{(2)}$,
    where $P^{(1)}$ is the subpath of $P$ which goes from the starting point of $P$ to the point $a$, and $P^{(2)}$ is the remaining subpath of $P$.
          Set $\beta(P)= \alpha(P^{(1)}) (P^{(2)})^{c}$.
 It is routine to check that   the resulting path $\beta(P)$ is a ballot $(s,k)$-path of   odd height.

 In order to show that the map $\beta$ is a bijection, we describe a map $\beta':\mathcal{B}^{o}(s,k) \rightarrow \mathcal{FB}_{-1}(s,k)$. Let $Q $ be a ballot $(s,k)$-path  of odd height $2\ell+1$ for some nonnegative  integer $\ell$ . Now we proceed to demonstrate how to get a path $\beta'(Q)$.  Let $b$ the rightmost lattice point of $P$ which is on the line $y=2\ell$. Then $Q$ can be uniquely decomposed as $Q=Q^{(1)}Q^{(2)}$, where $Q^{(1)}$ is the subpath of $Q$ which goes from the starting point of $Q$ to the point $b$, and $Q^{(2)}$ is the remaining subpath of $Q$. Set $\beta'(Q)= \alpha^{-1}(Q^{(1)})(Q^{(2)})^c $. One can easily verify that the resulting path $\beta'(Q)\in \mathcal{FB}_{-1}(s,k)$. Moreover, the ending point of $\alpha^{-1}(Q^{(1)})$ becomes the  rightmost lattice point which lies on the $x$-axis.    This ensures that the maps $\beta$ and $\beta'$ are inverses of each other. Thus the map $\beta$ is a bijection, completing the proof. \qed
\begin{example}
The decomposition of a path $$P=D_2U_2H_1U_2D_2H_2D_2H_1D_2H_3U_2H_1H_2$$ is shown in Figure \ref{beta1}. By applying the map $\beta$, its corresponding path is illustrated in Figure \ref{beta2}.
\end{example}

\begin{example}
The decomposition of a path $$Q=U_2H_1U_2D_2H_2U_2U_2H_1U_2H_3D_2H_1H_2$$ is shown in Figure \ref{beta2}. By applying the map $\beta'$, its corresponding path is illustrated in Figure \ref{beta1}.
\end{example}

 With abuse of notation, we use  $H_0$ to denote the empty step. Now we are in the position to establish a bijection between the set $\mathcal{Q}(m,k)$ and the set  $\mathcal{SD}(2m,2k)$.
\begin{thm}\label{delta}
There is a bijection between the  set $\mathcal{Q}(m,k)$ and the set  $\mathcal{SD}(2m,2k)$.
\end{thm}
\pf We first describe a map $\delta:  \mathcal{Q}(m, k) \rightarrow \mathcal{SD}(2m,2k) $. Let $P$ be a path in $\mathcal{Q}(m, k)$,  and  let $Q$ be the path obtained form $P$ by removing its first step.  We proceed to  construct a path $\delta(P)$ as follows.
\begin{itemize}
  \item[{\upshape (\rmnum{1})}] If $P\in \mathcal{FB}'_0(m-i,2k) $ for some $0\leq i\leq k-2$ and $P$ starts with a horizontal step $H_{\ell}$ for some $1\leq \ell\leq  2k-1$,  set $\delta(P)=\alpha(Q)H_{\ell}H_{2i}H_{\ell}  (\alpha(Q))^{rc} $.
  \item[{\upshape (\rmnum{2})}] If $P\in \mathcal{FB}'_0(m-i,2k) $ for some $0\leq i\leq k-2$ and $P$ starts with a down step $D_k$,  set $\delta(P)=\beta( Q ^{c})D_kH_{2i}U_k  (\beta( Q^{c}))^{rc}$.
     \item[{\upshape (\rmnum{3})}] If $P\in \mathcal{FB}'_{-1}(m+1,2k) $  and $P$ starts with a horizontal step $H_{\ell}$ for some $1\leq \ell\leq  k$,  set $\delta(P)=\beta(Q)H_{2\ell-2}  (\beta(Q))^{rc}$.
       \item[{\upshape (\rmnum{4})}] If $P\in \mathcal{FB}'_{-1}(m+1,2k) $  and $P$ starts with a horizontal step $H_{\ell}$ for some $k+1\leq \ell\leq 2k-1$,  set $\delta(P)=\beta(Q)U_kH_{2( \ell-k-1)}D_k(\beta(Q))^{rc}$.
           \item[{\upshape (\rmnum{5})}] If $P\in \mathcal{FB}'_{-1}(m+1,2k) $  and $P$ starts with a  down step $D_k$,  set $\delta(P)=\alpha(Q) H_{2k-2} (\alpha(Q))^{rc}$.
\end{itemize}
It is easy to check that  $\delta(P)\in SD(2m,2k)$, and thus the map $\delta$ is well defined.

In order to show the map $\delta$ is a bijection, we describe a map $\delta':  \mathcal{SD}(2m,2k) \rightarrow  \mathcal{Q}(m, k)$. Let $P$ be a path in $ \mathcal{SD}(2m,2k)$, we proceed to construct a path $\delta'(P)$ as follows.

\begin{itemize}
  \item[{\upshape ( \rmnum{1 }$'$ )}] If  $P= Q H_{\ell}H_{2i}H_{\ell}Q^{rc}$  with  $1\leq \ell\leq 2k-1$ and $0\leq i\leq k-2$ and   $Q$ is of even height, set $\delta'(P)=H_{\ell}\alpha^{-1}(Q)$.
  \item[{\upshape ( \rmnum{2 }$'$ )}] If  $P= Q  D_kH_{2i}U_k Q^{rc}$  with   $0\leq i\leq k-2$ and  $Q$ is of odd height, set $\delta'(P)=D_k(\beta^{-1}(Q))^c$.
     \item[{\upshape ( \rmnum{3 }$'$ )}] If    $P= Q   H_{2i} Q^{rc}$ with    $0\leq i\leq k-1$  and  $Q$ is of odd height, set $\delta'(P)=H_{i+1}\beta^{-1}(Q)$.
       \item[{\upshape ( \rmnum{4 }$'$ )}] If $P= Q  U_k H_{2i} D_k Q^{rc}$  with    $0\leq i\leq k-2$  and  $Q$ is of odd height, set $\delta'(P)=H_{k+1+i}\beta^{-1}(Q)$.
           \item[{\upshape ( \rmnum{5 }$'$ )}] If $P= Q   H_{2k-2}  Q^{rc}$  and  $Q$ is of even  height, set $\delta'(P)=D_k\alpha^{-1}(Q)$.
\end{itemize}
It is not difficult to check that the resulting path $\delta'(P)\in \mathcal{Q}(m,k)$. This implies that the map $\delta$ is well defined.  One can easily verify that {\upshape ( \rmnum{1 }$'$ )}, {\upshape ( \rmnum{2 }$'$ )},  {\upshape ( \rmnum{3 }$'$ )}, {\upshape ( \rmnum{4 }$'$ )} and {\upshape ( \rmnum{5 }$'$ )} respectively reverse the
procedures of {\upshape ( \rmnum{1 }  )}, {\upshape ( \rmnum{2 }  )},  {\upshape ( \rmnum{3 }  )}, {\upshape ( \rmnum{4 }  )} and {\upshape ( \rmnum{5 }  )}. This implies that the maps  $\delta$ and $\delta'$  are inverses of
each other. Hence, the map $\delta$ is a bijection, completing the proof. \qed

\begin{lem}\label{lemeven}
There is a bijection between $\mathcal{SD}(2m, 2k)$ and $\mathcal{SD}(2m+1,2k)$.
\end{lem}
 \pf  First we discribe a map $\xi:\mathcal{SD}(2m, 2k) \rightarrow \mathcal{SD}(2m+1,2k)$. For any path $P\in \mathcal{SD}(2m, 2k) $, it can be uniquely decomposed as $P=QH_{2\ell}Q^{rc}$ for some $0\leq \ell<{k}$. Let   $\xi(P)=QH_{2\ell+1}Q^{rc}$.  It is easily seen that the map $\xi$ is reversible, and thus is a bijection. This completes the proof.  \qed

  Theorems \ref{main1},  \ref{chi}, \ref{mainth2},  \ref{mainth3} and \ref{delta}  together with Lemma \ref{lemeven}  lead to   the following result.
\begin{thm}\label{theven}
For given positive integers $s$ and $ k$, the number of self-conjugate $(s, s+
1, \ldots, s + 2k)$-core partitions is equal to the number of symmetric $(s, 2k)$-Dyck paths.
\end{thm}

\subsection{On self-conjugate $(s, s+1, \ldots, s+2k+1)$-core partitions}
Let $$\mathcal{Q'}(m,k)=(\bigcup_{i=0}^{k-1}\mathcal{FB}'_0(m-i,2k+1))\cup \mathcal{FB}'_{-1}(m+{1\over 2},2k+1).$$
In the plane,   a point $(x,y)$  is said to be {\em visible} if  either $x=a+{1\over 2}$ for some $a\in \mathcal{Z}$ and $y$ is an  odd integer, or $x=a$ for some $a\in \mathcal{Z}$ and $y$ is an  even integer.

 We describe a map $\psi:J(P'(2m,2k+1))\rightarrow \mathcal{Q'}(m,k)$ as follows.  Let $I\in J(P'(2m,2k+1)) $. For the visible  point $(a,b)$  in the plane  where $ b\in \mathcal{Z}$ and   $0\leq a\leq m+{1\over 2}$, we color it  by black and white according to the following rules.
\begin{itemize}
\item If $b\geq 0$ , then
    the point $(a,b)$  is colored by  black if and only   if $(a,b) \in I  $;
    \item If $b<0$, then the point $(a,b)$  is colored by  black if and only    if $(a,b)\notin I  $.
    \end{itemize}

\begin{figure}[h!]
\begin{center}
\begin{picture}(280,160)
\setlength{\unitlength}{8mm}

 \put(0,0){\circle*{0.2}} \put(1,0){\circle*{0.2}} \put(2,0){\circle*{0.2}} \put(3,0){\circle*{0.2}} \put(4,0){\circle*{0.2}} \put(5,0){\circle*{0.2}}
  \put(6,0){\circle*{0.2}} \put(7,0){\circle*{0.2}} \put(8,0){\circle*{0.2}}  \put(9,0){\circle*{0.2}} \put(10,0){\circle*{0.2}}

  \put(0.5,1){\circle*{0.2}} \put(1.5,1){\circle*{0.2}} \put(2.5,1){\circle*{0.2}} \put(3.5,1){\circle*{0.2}} \put(4.5,1){\circle*{0.2}} \put(5.5,1){\circle*{0.2}}
  \put(6.5,1){\circle*{0.2}} \put(7.5,1){\circle*{0.2}} \put(8.5,1){\circle*{0.2}}  \put(9.5,1){\circle*{0.2}} \put(10.5,1){\circle*{0.2}}

\put(0,2){\circle*{0.2}} \put(1,2){\circle*{0.2}} \put(2,2){\circle*{0.2}} \put(3,2){\circle*{0.2}} \put(4,2){\circle*{0.2}} \put(5,2){\circle*{0.2}}
  \put(6,2){\circle{0.2}} \put(7,2){\circle*{0.2}} \put(8,2){\circle*{0.2}}  \put(9,2){\circle*{0.2}} \put(10,2){\circle*{0.2}}

  \put(0.5,3){\circle*{0.2}} \put(1.5,3){\circle*{0.2}} \put(2.5,3){\circle*{0.2}} \put(3.5,3){\circle*{0.2}} \put(4.5,3){\circle{0.2}} \put(5.5,3){\circle{0.2}}
  \put(6.5,3){\circle{0.2}} \put(7.5,3){\circle{0.2}} \put(8.5,3){\circle*{0.2}}  \put(9.5,3){\circle*{0.2}} \put(10.5,3){\circle*{0.2}}

  \put(0,4){\circle{0.2}} \put(1,4){\circle*{0.2}} \put(2,4){\circle{0.2}} \put(3,4){\circle{0.2}} \put(4,4){\circle{0.2}} \put(5,4){\circle{0.2}}
  \put(6,4){\circle{0.2}} \put(7,4){\circle{0.2}} \put(8,4){\circle{0.2}}  \put(9,4){\circle{0.2}} \put(10,4){\circle{0.2}}

  \put(0.5,5){\circle{0.2}} \put(1.5,5){\circle{0.2}} \put(2.5,5){\circle{0.2}} \put(3.5,5){\circle{0.2}} \put(4.5,5){\circle{0.2}} \put(5.5,5){\circle{0.2}}
  \put(6.5,5){\circle{0.2}} \put(7.5,5){\circle{0.2}} \put(8.5,5){\circle{0.2}}  \put(9.5,5){\circle{0.2}} \put(10.5,5){\circle{0.2}}

  \put(0,6){\circle{0.2}} \put(1,6){\circle{0.2}} \put(2,6){\circle{0.2}} \put(3,6){\circle{0.2}} \put(4,6){\circle{0.2}} \put(5,6){\circle{0.2}}
  \put(6,6){\circle{0.2}} \put(7,6){\circle{0.2}} \put(8,6){\circle{0.2}}  \put(9,6){\circle{0.2}} \put(10,6){\circle{0.2}}

\put(0,4){\line(1,0){1}}\put(1,4){\line(3,-2){1.5}}\put(2.5,3){\line(1,0){1}}\put(3.5,3){\line(3,-2){1.5}}\put(5,2){\line(2,0){2}} \put(7,2){\line(3,2){1.5}}\put(8.5,3){\line(1,0){1}}\put(9.5,3){\line(1,0){1}}

\put(-1.2,3.8){\small$(0,0)$}

\end{picture}
\end{center}
\caption{ The corresponding path $\psi(I)$.}\label{psi}
\end{figure}
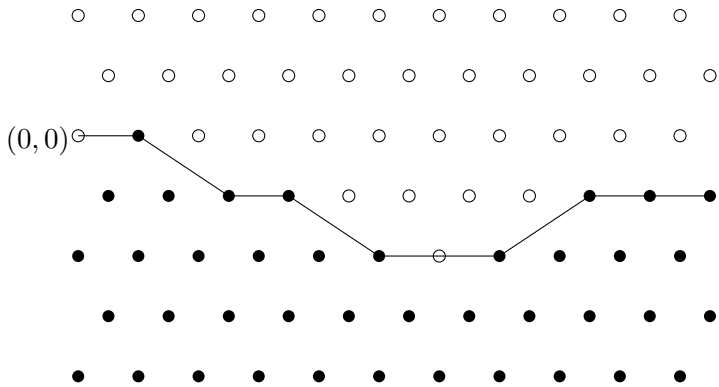

 Now we recursively  define a sequence   $P_0,P_1, \ldots P_n$ of points where  $P_0=(0,0)$. Assume that $  P_{i-1}$ has  been already determined. Let   $P_{i-1}=(a_{i-1}, b_{i-1})$.   We proceed to  get $P_i$ by the following procedure.\\
 Case 1. There exists a black point $(a_{i-1}+\ell, b_{i-1})$  with    $1\leq \ell\leq  2k$.  If the point $ (a_{i-1}+k+{1\over 2}, b_{i-1}+1)$ is black, then  set $P_i=(a_{i-1}+k+{1\over 2}, b_{i-1}+1)$. Otherwise,  choose $\ell$ to the smallest such integer and set  $P_i=(a_{i-1}+\ell, b_{i-1})$.\\
 Case 2. There does not exist  a black point $(a_{i-1}+\ell, b_{i-1})$  with    $1\leq \ell\leq  2k$. If the point $ (a_{i-1}+k+{1\over 2}, b_{i-1}-1)$ is black, set $P_i=(a_{i-1}+k+{1\over 2}, b_{i-1}-1)$. Otherwise, set $P_{i}=P_{i-1}$ and  $n=i-1$.\\
Let $\psi(I)$ be the resulting  path with lattice points $P_0, P_1, \ldots, P_n$.
For example, let $I=\{(1,0), (6,-2), (7.5,-1), (6.5, -1), (5.5,-1), (4.5,-1)\}\in J(P'(20,3))$. By applying the map $\psi$, we get a  path $\psi(I)$ as shown in Figure \ref{psi}.

By the same reasoning  as in the proof of Lemmas \ref{lemphi3}, \ref{lemphi1}, \ref{lemphi2} and Theorem \ref{mainth3}, we can deduce the following results.

   \begin{lem}\label{lempsi1}
 The map $\psi$ is well defined, that is, for any $I\in J(P'(2m,2k+1))$  we have   $ \psi(I)\in\mathcal{Q}'(m,k)$.
\end{lem}

\begin{lem}\label{lempsi2}
For any   $I\in J(P'(2m,2k+1))$,  the resulting path $\psi(I)$ verifies the following properties.
\begin{itemize}
\item[{\upshape (\rmnum{1})}] The visible points   on  the path $\psi(I)$  are black if and only if they are ending points of steps ;
\item[ {\upshape (\rmnum{2})}] For $b<0$, the visible point  $(a,b)$       to the right of the ending point of  the path $\psi(I)$ is black ;
\item[ {\upshape (\rmnum{3})}] For $b\geq 0$, the visible point $(a,b)$        to the right of the ending point of  the path $\psi(I)$ is white;
 \item[{\upshape (\rmnum{4})}] The visible points   lying     above the path $\psi(I)$ are white;
\item[ {\upshape (\rmnum{5})}] The visible points   lying     below the path $\psi(I)$ are black;

\end{itemize}
\end{lem}

\begin{thm}\label{mainth4}
The map $\psi$ is a bijection between  $J(P'(2m,2k+1)) $ and $\mathcal{Q}'(m,k)$.
\end{thm}

In the following theorem, we establish a bijection between  the set $\mathcal{Q}'(m,k)$ and the set  $\mathcal{SD}(2m,2k+1)$.

\begin{thm}\label{gamma}
There is a bijection   between the  set $\mathcal{Q}'(m,k)$ and the set  $\mathcal{SD}(2m,2k+1)$.
\end{thm}
\pf We first describe a map $\gamma:  \mathcal{Q}'(m,k) \rightarrow \mathcal{SD}(2m,2k+1) $. Let $P$ be a path in $\mathcal{Q}'(m,k) $ and let $Q$ be the path obtained from $P$ by removing its first step. Now  we proceed to demonstrate the procedure to get a path $\gamma(P)$.
\begin{itemize}
  \item[{\upshape (\rmnum{1})}] If $P\in \mathcal{FB}'_0(m-i,2k+1))$ for some $0\leq i\leq k-1$ and $P$ starts with a horizontal step $H_{\ell}$ for some $1\leq \ell\leq  2k$,  set $\gamma(P)=\alpha(Q)H_{\ell}H_{2i}H_{\ell}(\alpha(Q))^{rc}$.
  \item[{\upshape (\rmnum{2})}] If $P\in \mathcal{FB}'_0(m-i,2k+1))$ for some $0\leq i\leq k-1$ and $P$ starts with a down step $D_{k+{1\over 2}}$,  set $\gamma(P)=\beta(Q^c)D_{k+{1\over 2}}H_{2i}U_{k+{1\over 2}}(\beta(Q^c))^{rc}$.
     \item[{\upshape (\rmnum{3})}] If $P\in \mathcal{FB}'_{-1}(m+{1\over 2},2k+1))$  and $P$ starts with a horizontal step $H_{\ell}$ for some $1\leq \ell\leq  k$,  set $\gamma(P)=\beta(Q)H_{2\ell-1} (\beta(Q))^{rc}$.
       \item[{\upshape (\rmnum{4})}] If $P\in \mathcal{FB}'_{-1}(m+{1\over 2},2k+1))$  and $P$ starts with a horizontal step $H_{\ell}$ for some $k+1\leq \ell\leq 2k$,  set $\gamma(P)=\beta(Q)U_{k+{1\over 2}}H_{2( \ell-k-1)}D_{k+{1\over 2}}(\beta(Q))^{rc}$.
           \item[{\upshape (\rmnum{5})}] If $P\in \mathcal{FB}'_{-1}(m+{1\over 2},2k+1))$  and $P$ starts with a  down step $D_{k+{1\over 2}}$,  set $\gamma(P)=\alpha(Q) H_{2k} (\alpha(Q))^{rc} $.
\end{itemize}
One can easily  check that  $\gamma(P)\in \mathcal{SD}(2m,2k+1)$, and thus the map $\gamma$ is well defined.

Next we describe a map $\gamma':  \mathcal{SD}(2m,2k+1) \rightarrow  \mathcal{Q}'(m,k) $. Let $P$ be a path in $ \mathcal{SD}(2m,2k+1)$, we proceed to construct a path $\gamma'(P)$ by the following procedure.

\begin{itemize}
  \item[{\upshape ( \rmnum{1 }$'$ )}] If  $P= Q H_{\ell}H_{2i}H_{\ell}Q^{rc}$ with  $1\leq \ell\leq 2k$ and $0\leq i\leq k-1$ and $Q$ is of even height, set $\gamma'(P)=H_{\ell}\alpha^{-1}(Q)$.
  \item[{\upshape ( \rmnum{2 }$'$ )}] If  $P= Q  D_{k+{1\over 2}}H_{2i}U_{k+{1\over 2}}Q^{rc}$ with    $0\leq i\leq k-1$  and $Q$ is of odd height, set $\gamma'(P)=D_{k+{1\over 2}}(\beta^{-1}(Q))^c$.
     \item[{\upshape ( \rmnum{3 }$'$ )}] If    $P= Q   H_{2i-1} Q^{rc}$ with   $1\leq i\leq k$  and  $Q$ is of odd height, set $\gamma'(P)=H_{i}\beta^{-1}(Q)$.
       \item[{\upshape ( \rmnum{4 }$'$ )}] If $P= Q  U_{k+{1\over 2}} H_{2i} D_{k+{1\over 2}} Q^{rc}$  with  $0\leq i\leq k-1$  and  $Q$ is of odd height, set $\gamma'(P)=H_{k+1+i}\beta^{-1}(Q)$.
           \item[{\upshape ( \rmnum{5 }$'$ )}] If $P= Q    H_{2k}  Q^{rc}$  and  $Q$ is of even  height, set $\gamma'(P)=D_{k+{1\over 2}}\alpha^{-1}(Q)$.
\end{itemize}
It is not difficult to check that the resulting path $\gamma'(P)\in \mathcal{Q}'(m,k) $. This implies that the map $\gamma$ is well defined. Furthermore, one can easily verify that {\upshape ( \rmnum{1 }$'$ )}, {\upshape ( \rmnum{2 }$'$ )},  {\upshape ( \rmnum{3 }$'$ )}, {\upshape ( \rmnum{4 }$'$ )} and {\upshape ( \rmnum{5 }$'$ )} respectively reverse the
procedures of {\upshape ( \rmnum{1 }  )}, {\upshape ( \rmnum{2 }  )},  {\upshape ( \rmnum{3 }  )}, {\upshape ( \rmnum{4 }  )} and {\upshape ( \rmnum{5 }  )}. This implies that the maps  $\gamma$ and $\gamma'$  are inverses of
each other. Hence, the map $\gamma$ is a bijection, completing the proof. \qed

Combining   Theorems \ref{main1},   \ref{chi}, \ref{mainth4} and \ref{gamma}, we are    led  to   the following result.
\begin{thm}\label{thodd}
For given positive integer  $m$ and  nonnegative integer $ k$, the number of self-conjugate $(2m, 2m+
1, \ldots, 2m + 2k+1 )$-core partitions is equal to the number of symmetric $(2m, 2k+1)$-Dyck paths.
\end{thm}

Let $ J^*(P'(2m+2, 2k+1))$ denote   the order ideals $ I \in  J(P'(2m+2, 2k+1))$ such that $I$ does not contain the element $(m+1, 0)$.   It is easily seen that \begin{equation}\label{eq1}
J^*(P'(2m+2, 2k+1))= J(P'(2m+1, 2k+1)).
\end{equation}
Let  $\mathcal{Q}^*(m,k)=\mathcal{Q}'(m,k)\setminus \mathcal{FB}'_0(m, 2k+1)$. The construction of the map $\psi$ together with Lemma \ref{lempsi2} leads to  the following result.

\begin{thm}\label{psi'}
The map $\psi$ induces a bijection between the set $J^*(P'(2m+2, 2k+1))$ and the set  $\mathcal{Q}^*(m+1,k)$.
\end{thm}

In the following theorem, we aim to establish a bijection between the  set $\mathcal{Q}^*(m+1,k)$ and the set  $\mathcal{SD}(2m+1,2k+1)$.
\begin{thm}\label{eta}
There is a bijection   between the  set $\mathcal{Q}^*(m+1,k)$ and the set  $\mathcal{SD}(2m+1,2k+1)$.
\end{thm}
\pf We first describe a map $\eta:  \mathcal{Q}^*(m+1,k) \rightarrow \mathcal{SD}(2m+1,2k+1) $. Let $P$ be a path in $\mathcal{Q}^*(m+1,k) $ and let $Q$ be the path obtained from $P$ by removing its first step. Now  we proceed to demonstrate the procedure to get a path $\eta(P)$.
\begin{itemize}
  \item[{\upshape (\rmnum{1})}] If $P\in \mathcal{FB}'_0(m+1-i,2k+1))$ for some $1\leq i\leq k-1$ and $P$ starts with a horizontal step $H_{\ell}$ for some $1\leq \ell\leq  2k$,  set $\eta(P)=\alpha(Q)H_{\ell}H_{2i-1}H_{\ell}(\alpha(Q))^{rc}$.
  \item[{\upshape (\rmnum{2})}] If $P\in \mathcal{FB}'_0(m+1-i,2k+1))$ for some $1\leq i\leq k-1$ and $P$ starts with a down step $D_{k+{1\over 2}}$,  set $\eta(P)=\beta(Q^c)D_{k+{1\over 2}}H_{2i-1}U_{k+{1\over 2}}(\beta(Q^c))^{rc}$.
     \item[{\upshape (\rmnum{3})}] If $P\in \mathcal{FB}'_{-1}(m+1+{1\over 2},2k+1))$  and $P$ starts with a horizontal step $H_{\ell}$ for some $1\leq \ell\leq  k+1$,  set $\eta(P)=\beta(Q)H_{2\ell-2} (\beta(Q))^{rc}$.
       \item[{\upshape (\rmnum{4})}] If $P\in \mathcal{FB}'_{-1}(m+1+{1\over 2},2k+1))$  and $P$ starts with a horizontal step $H_{\ell}$ for some $k+2\leq \ell\leq 2k$,  set $\eta(P)=\beta(Q)U_{k+{1\over 2}}H_{2( \ell-k-1)-1}D_{k+{1\over 2}}(\beta(Q))^{rc}$.
           \item[{\upshape (\rmnum{5})}] If $P\in \mathcal{FB}'_{-1}(m+1+{1\over 2},2k+1))$  and $P$ starts with a  down step $D_{k+{1\over 2}}$,  set $\eta(P)=\alpha(Q) H_{2k-1} (\alpha(Q))^{rc} $.
\end{itemize}
It is easy to verify that $\eta(P)\in \mathcal{SD}(2m+1,2k+1)$, and thus the map $\eta$ is well defined.

Next we describe a map $\eta':  \mathcal{SD}(2m+1,2k+1) \rightarrow  \mathcal{Q}^*(m+1,k) $. Let $P$ be a path in $ \mathcal{SD}(2m+1,2k+1)$, we proceed to construct a path $\eta'(P)$ by the following procedure.

\begin{itemize}
  \item[{\upshape ( \rmnum{1 }$'$ )}] If  $P= Q H_{\ell}H_{2i-1}H_{\ell}Q^{rc}$ with  $1\leq \ell\leq 2k$ and $1\leq i\leq k-1$ and $Q$ is of even height, set $\eta'(P)=H_{\ell} \alpha^{-1}(Q)$.
  \item[{\upshape ( \rmnum{2 }$'$ )}] If  $P= Q  D_{k+{1\over 2}}H_{2i-1}U_{k+{1\over 2}}Q^{rc}$ with    $1\leq i\leq k-1$  and $Q$ is of odd height, set $\eta'(P)=D_{k+{1\over 2}}(\beta^{-1}(Q))^c$.
     \item[{\upshape ( \rmnum{3 }$'$ )}] If    $P= Q   H_{2i} Q^{rc}$ with   $0\leq i\leq k$  and  $Q$ is of odd height, set $\eta'(P)=H_{i+1}\beta^{-1}(Q)$.
       \item[{\upshape ( \rmnum{4 }$'$ )}] If $P= Q  U_{k+{1\over 2}} H_{2i-1} D_{k+{1\over 2}} Q^{rc}$  with  $1\leq i\leq k-1$  and  $Q$ is of odd height, set $\eta'(P)=H_{k+1+i}\beta^{-1}(Q)$.
           \item[{\upshape ( \rmnum{5 }$'$ )}] If $P= Q    H_{2k-1}  Q^{rc}$  and  $Q$ is of even  height, set $\eta'(P)=D_{k+{1\over 2}}\alpha^{-1}(Q)$.
\end{itemize}
It is routine  to check that the resulting path $\eta'(P)\in \mathcal{Q}^*(m+1,k) $. Therefore, the map $\eta'$ is well defined. Furthermore, one can easily verify that {\upshape ( \rmnum{1 }$'$ )}, {\upshape ( \rmnum{2 }$'$ )},  {\upshape ( \rmnum{3 }$'$ )}, {\upshape ( \rmnum{4 }$'$ )} and {\upshape ( \rmnum{5 }$'$ )} respectively reverse the
procedures of {\upshape ( \rmnum{1 }  )}, {\upshape ( \rmnum{2 }  )},  {\upshape ( \rmnum{3 }  )}, {\upshape ( \rmnum{4 }  )} and {\upshape ( \rmnum{5 }  )}. This implies that the maps  $\eta$ and $\eta'$  are inverses of
each other. Hence, the map $\eta$ is a bijection, completing the proof. \qed

 Theorems \ref{main1},   \ref{chi}, \ref{psi'} and \ref{eta} together with Equation (\ref{eq1})  yield  the following result.
\begin{thm}\label{thodd1}
For given positive integer  $m$ and  nonnegative integer $ k$, the number of self-conjugate $(2m+1, 2m+
2, \ldots, 2m + 2k+2 )$-core partitions is equal to the number of symmetric $(2m+1, 2k+1)$-Dyck paths.
\end{thm}

 Combining Theorems \ref{theven},  \ref{thodd} and \ref{thodd1}, we complete the  proof of Conjecture \ref{con1}.

\noindent{\bf Acknowledgments.} This work was supported by  the National Natural Science Foundation of China (11671366).

%--------------------------

\end{document}